\documentclass[letterpaper]{article} 
\usepackage{arxiv}  
\usepackage{times}  
\usepackage{helvet}  
\usepackage{courier}  
\usepackage[hyphens]{url}  
\usepackage{graphicx} 
\usepackage{caption} 
\usepackage[utf8]{inputenc} 
\usepackage[T1]{fontenc}    
\usepackage{booktabs}       
\usepackage{amsfonts}       
\usepackage{nicefrac}       
\usepackage{microtype}      
\usepackage{lipsum}         
\usepackage{subcaption}
\usepackage{amsmath}
\usepackage{algorithm}
\usepackage{algpseudocode}
\usepackage{pgfgantt}
\usepackage{appendix}
\usepackage{mathrsfs}
\usepackage{multirow}
\newtheorem{remark}{Remark}
\usepackage{adjustbox}
\usepackage{threeparttable}
\usepackage{hyperref}       
\usepackage{tabularx}

\makeatletter
\newsavebox{\@tabnotebox}
\providecommand\tmark{} 
\providecommand\tnote{}
\newenvironment{tabularwithnotes}[3][c]
{\long\def\@tabnotes{#3}%
\renewcommand\tmark[1][a]{\makebox[0pt][l]{\textsuperscript{\itshape##1}}}%
\renewcommand\tnote[2][a]{\textsuperscript{\itshape##1}\,##2\par}
\begin{lrbox}{\@tabnotebox}
\begin{tabular}{#2}}
{\end{tabular}\end{lrbox}%
\parbox{\wd\@tabnotebox}{
\usebox{\@tabnotebox}\par
\smallskip\@tabnotes
}%
}
\makeatother
\title{Solving Integrated Process Planning and Scheduling Problem via Graph Neural Network Based Deep Reinforcement Learning}

\author{
    \and Hongpei Li$^{1,*}$\and  Han Zhang$^{1,*}$\and  Ziyan He$^1$\and  Yunkai Jia$^2$\and Bo Jiang$^1$\and  Xiang Huang$^2$\and  Dongdong Ge$^3$
    \\
    \\
    \textsuperscript{1}Shanghai University of Finance and Economics\\
    \textsuperscript{2}Cardinal Operations\\
    \textsuperscript{3}Shanghai Jiao Tong University
}

\begin{document}
\maketitle
\footnotetext[1]{These authors contributed equally to this work.}
\footnotetext[2]{Correspondence to: Hongpei Li<ishongpeili@gamil.com>, Bojiang<isyebojiang@gmail.com>}
\begin{abstract}
    The Integrated Process Planning and Scheduling (IPPS) problem combines process route planning and shop scheduling to achieve high efficiency in manufacturing and maximize resource utilization, which is crucial for modern manufacturing systems. Traditional methods using Mixed Integer Linear Programming (MILP) and heuristic algorithms can not well balance solution quality and speed when solving IPPS. In this paper, we propose a novel end-to-end Deep Reinforcement Learning (DRL) method. We model the IPPS problem as a Markov Decision Process (MDP) and employ a Heterogeneous Graph Neural Network (GNN) to capture the complex relationships among operations, machines, and jobs. To optimize the scheduling strategy, we use Proximal Policy Optimization (PPO). Experimental results show that, compared to traditional methods, our approach significantly improves solution efficiency and quality in large-scale IPPS instances, providing superior scheduling strategies for modern intelligent manufacturing systems.
\end{abstract}


\section{Introduction}
Deep Reinforcement Learning (DRL) is a popular technique to learn sequential decision-making in 
complex problems, combining reinforcement learning with deep learning \cite{SHAKYA2023120495}, which has achieved successful applications in various areas including robotics, traffic signal control, and recommendation systems. There is a growing interest in applying DRL to NP-hard combinatorial optimization problems (COPs) \cite{cappart2021combining} to find a reasonable good solution within an acceptable time such as the Traveling Salesman Problem (TSP) and scheduling problems such as Job Shop Scheduling (JSP) and  Flexible Job Shop Scheduling Problem (FJSP) \cite{zhang2020deep,Zhang2023,Liu2023}. Since DRL inherits the fundamental framework of dynamic programming and Markov Decision Process (MDP), it is naturally suitable for solving multi-stage decision problems transformed from COPs. 

Inspired by the successful applications mentioned above, we introduce DRL to the Integrated Process Planning and Scheduling (IPPS) problem, a typical NP-hard COP that is crucial for intelligent manufacturing. IPPS combines two of the most important subsystems in manufacturing systems: process planning and shop scheduling, where the former plans the optimal process route under precedence constraints, and the latter schedules machine assignments, sequences of jobs, and start times based on the known process route. Traditionally, these two subsystems execute sequentially. However, performing them separately can lead to bottleneck resources, conflicts in optimization objectives, and unbalanced machine loads \cite{Liu2021}. In contrast, IPPS simultaneously optimizes both process routes
and scheduling plans under specific constraints and objectives, such as completion time \cite{Phanden2021}. Therefore, IPPS has become a popular topic in current production research.

IPPS has been widely applied in the automotive, shipbuilding, chemical, metallurgy, and aerospace manufacturing areas \cite{doi:10.1080/00207543.2021.1892230,samaranayake2012aircraft,varela2021integrated} which require modern manufacturing systems with efficiently and dynamically responding to rapidly changing product demands\cite{Zhou2019}. However, as 
a typical NP-hard problem, traditional methods for solving IPPS, such as mixed-integer linear programming and heuristic genetic algorithms \cite{jin2016,barzanji2020,kim2003,li2018}, struggle to quickly adapt to dynamic, real-world environments with unexpected disruptions due to their computational complexity. Therefore, ensuring the efficiency of IPPS solutions has become a challenging issue.

In combinatorial optimization problems similar to IPPS, DRL methods show advantages in efficiency and accuracy. Taking the FJSP as an example, DRL models the scheduling process as an MDP, using neural network models to collect production environment information and output scheduling priorities. \cite{luo2021} use Priority Dispatching Rules (PDRs) derived from expert experience and employ deep Q-networks to solve FJSP end-to-end. To fully exploit the graph structure of FJSP, \cite{song2023,echeverria2023} model the state of FJSP as a heterogeneous graph, combining graph neural networks with attention mechanisms, which perform well on large-scale problems. Unlike heuristic methods, DRL has strong transfer learning capabilities, quickly adapts to new problems, and provides high-quality solutions.

However, current research rarely uses DRL to solve IPPS. The existing DRL methods used for FJSP cannot apply to IPPS directly due to additional complexity in the decision-making and representation of IPPS while simultaneously solving process route selection and resource scheduling. Thus, in this paper, we propose a novel end-to-end DRL method. We model the IPPS problem as an MDP and employ a Heterogeneous Graph Neural Network (HGNN) to capture the complex relationships among operations, machines, and jobs. To optimize the scheduling strategy, we use Proximal Policy Optimization (PPO).

\textbf{Contributions}  To the best of our knowledge, this article is the first to apply DRL to solve the IPPS problem. We formulate an equivalent MDP for IPPS by incorporating an empty action to represent waiting at decision points. Presenting the current state as a heterogeneous graph, we use GNN for feature extraction and generate an operation-machine pair to process next by DRL structure.
\begin{itemize}
    \item \textbf{Exploration Efficiency:} We introduce the concept of combination—a set of indispensable operations for finishing the job—into our action space reduction strategy. While the number of combinations is relatively small, ineffective exploration costs associated with multiple production paths are avoided.
    \item \textbf{Stability and Generalization:} Our state transition process prunes completed or infeasible operations and their connections from the heterogeneous graph, reducing noise and focusing on critical scheduling tasks. Moreover, we propose a dense reward function based on estimated completion times, which outperforms a naive sparse reward related to the actual completion time. Together, our model can generalize effectively across various scheduling environments through a stable training process.
    \item \textbf{Competitive Numerical Performance: }Training on small-size problems with no more than 6 jobs, our DRL approach substantially outperforms in large-scale problems, particularly those with more than 16 jobs over 30,000 variables, achieving an improvement of 11.35\% compared to OR-Tools.

\end{itemize}

\section{Related Work}
\subsection{Conventional IPPS Methods}

Traditional IPPS solutions include exact methods like MILP and Constraint Programming (CP). \cite{jin2016} developed MILP models for IPPS using network representations, solved with CPLEX. \cite{barzanji2020} used logic-based benders decomposition to iteratively solve sub-problems, but these approaches struggle with large-scale problems. Consequently, research has shifted towards heuristic methods.

\cite{kim2003} establish a benchmark dataset using a symbiotic evolutionary algorithm. \cite{li2007} introduce simulated annealing for optimization, improving performance on Kim's dataset. Other hybrid methods, such as ant colony optimization \cite{leung2010}, tabu search \cite{li2010}, variable neighborhood search \cite{li2018}, and harmony search \cite{wu2021}, have been applied to IPPS. While faster than exact algorithms, these approaches lack optimality guarantees.

\subsection{DRL-based Scheduling Methods}

To enhance solution efficiency and quality, DRL has been applied to scheduling problems like FJSP. \cite{luo2021} combined expert-based Priority Dispatching Rules (PDR) with deep Q-Networks for dynamic FJSP, while \cite{han2021} proposed an end-to-end DRL method using tripartite graphs. However, these approaches overlook graph structures in scheduling problems.

Recent studies incorporate Graph Neural Networks (GNNs) to capture unique graph information in FJSP. \cite{zhang2020} used interlaced graphs and DRL to learn high-quality PDR, while \cite{song2023} modeled FJSP as an MDP with heterogeneous graphs and attention mechanisms for state representation. \cite{echeverria2023} utilized multi-strategy generation and PPO with PDR to restrict the action space, achieving strong performance in large-scale instances.

Most related works focus on using DRL for FJSP, but these methods can't be directly applied to IPPS. Improving solving IPPS with DRL remains a crucial research topic.

\section{Preliminaries}
\subsection{Graph-based IPPS Problem}

An IPPS instance is defined by a set of \( n \) jobs \( \mathcal{J} = \{ J_1, J_2, \ldots, J_n \} \) and a set of \( m \) machines \( \mathcal{M} = \{ M_1, M_2, \ldots, M_m \} \). For each job \( J_i \in \mathcal{J} \), \( J_i \) has a set of operations \( \mathcal{O}_i = \{ O_{i1}, O_{i2}, \ldots, O_{il} \} \), and each operation \( O_{ij} \) in \( \mathcal{O}_i \) can be processed on machines belonging to a set of machines \( \mathcal{M}_{ij} \subseteq \mathcal{M} \). The time machine \( M_k \) processes the operation \( O_{ij} \) is denoted as \( p_{ijk} \in \mathbb{R}^+ \). The most common objective function used to minimize an IPPS problem is makespan, which is defined as \(\max_{i,j} T_{ij}\) where \(T_{ij}\) is the completion time of \(O_{ij}\).

\begin{figure}[ht]
\centering
    \begin{subfigure}[b]{0.24\textwidth}
        \centering
        \includegraphics[width=\textwidth]{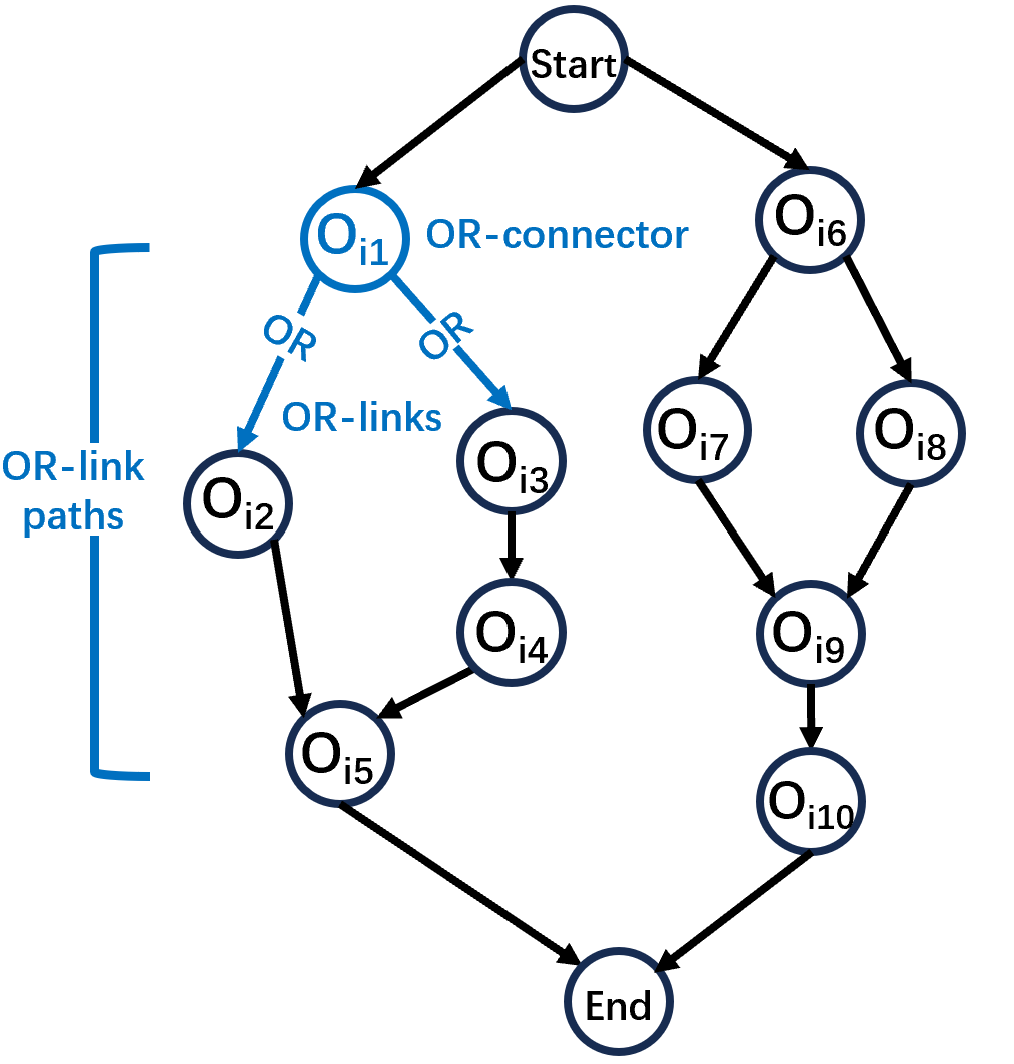}
        \caption{}
        \label{fig:andor}
    \end{subfigure}
    \begin{subfigure}[b]{0.20\textwidth}

        \includegraphics[width=\textwidth]{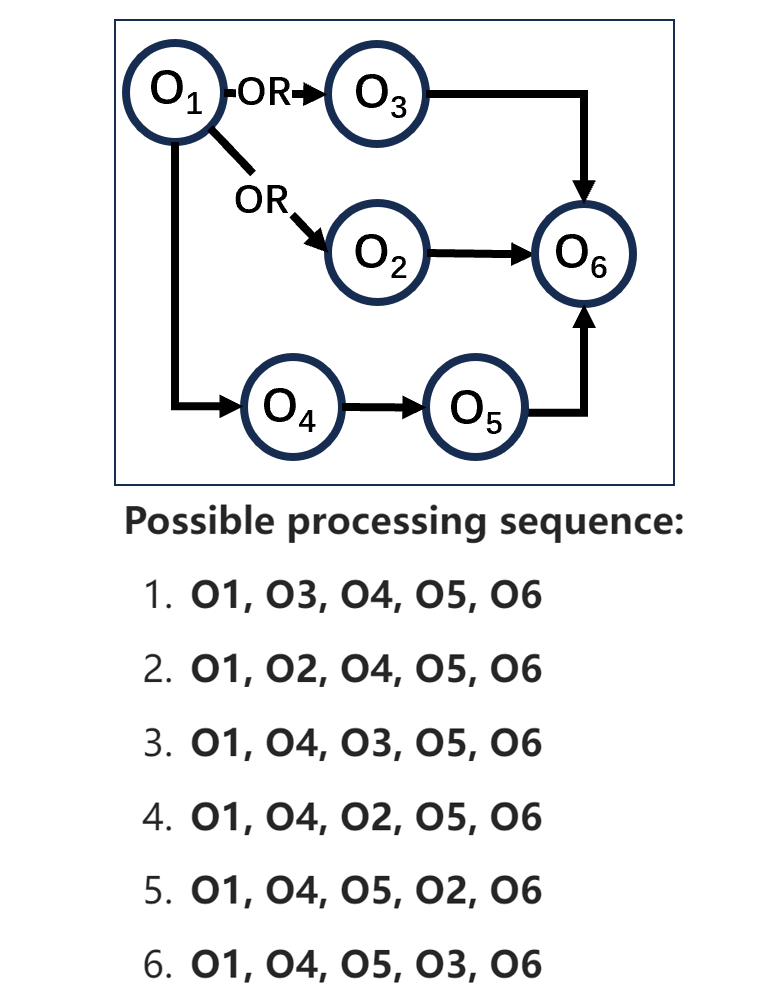}
        \caption{}
        \label{fig:seqs}
    \end{subfigure}
    \caption{(a) AND/OR graph of an IPPS job. (b) All possible processing sequences for a single instance.}
    \label{fig:overall}
\end{figure}

To complete each job \( J_i \in \mathcal{J} \), operations in \( \mathcal{O}_i  \) must follow certain precedence constraints, which can be represented in an AND/OR graph \cite{ho1996solving}. In an AND/OR graph, each job begins and ends with a starting node and an ending node, and intermediate nodes represent each operation in \( \mathcal{O}_i \)(see Figure~\ref{fig:andor}). The arrows connecting the nodes represent the precedence between them. Arrows with the attribute 'OR' are called OR-links, and the start nodes of OR-links are OR-connectors. An operation path that begins at an OR-link and ends when it merges with other paths is called an OR-link path. When processing jobs, only one path in OR-paths beginning with the same OR-connectors is chosen. For nodes that are not in any OR-path, all of them should be visited. Following the rules above, we define the set of indispensable operations for finishing the job without considering precedence relationships between operations as a combination. \( \mathcal{C}_i = \bigcup\limits_{h} \mathcal{C}_{ih} \) is the set of all possible combinations of \( J_i \).
\begin{remark}
The definition of a combination is proposed by \cite{jin2016}. We mainly focus on combination without precedence relationships rather than each processing sequence which shows a specific processing order in this paper, because the number of combinations is much less than the processing sequences that can be up to \(O(n!)\), and provides an effective exploration space.
\end{remark}

\section{Motivation and Novel Ingredient}
\textbf{Motivation: }In fact, DRL has been utilized to solve some easier scheduling problems like JSP and FJSP, which however can not be adapted to IPPS. In particular, each job in FJSP has a known single processing sequence, but the processing sequence for each job in IPPS could have \(O(n!)\) permutations (see Figure~\ref{fig:seqs}). If we consider all possible sequences in IPPS, then solving one IPPS can be reduced to solving \(O(n!)\) FJSP, which becomes computationally infeasible even if DRL for FJSP is fast. Thus directly applying previous DRL methods of FJSP to IPPS is impractical.

\begin{figure}[ht]
    \centering
    \includegraphics[width=0.5\linewidth]{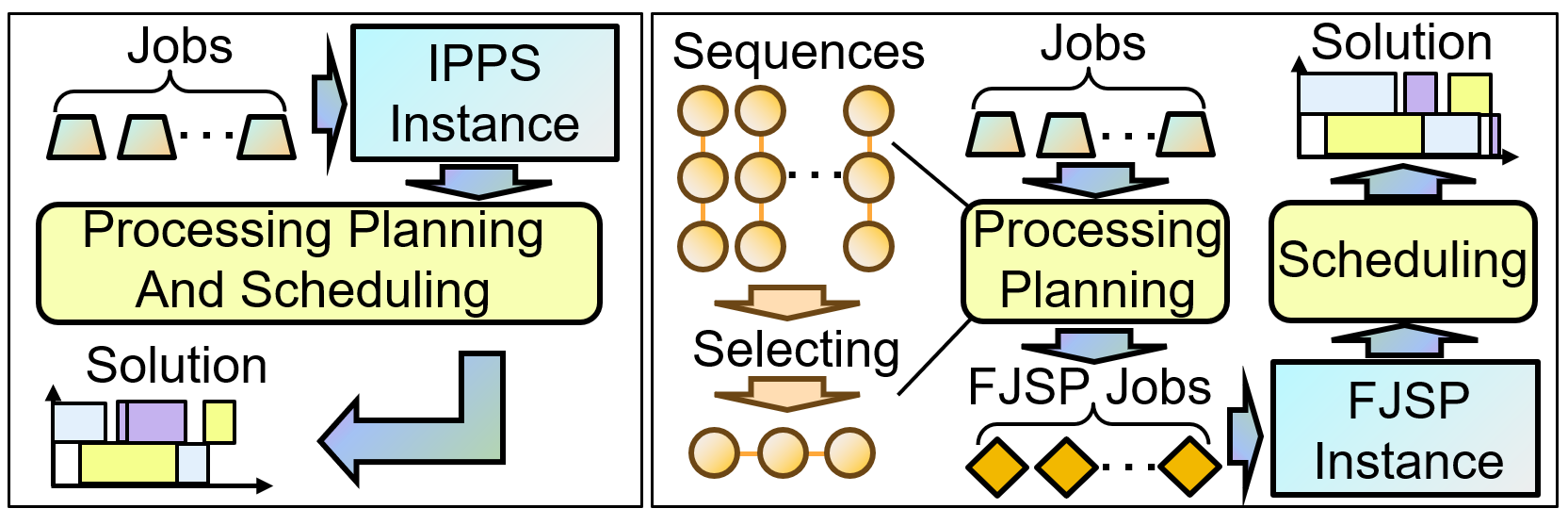}
    \caption{Comparison between IPPS (left) and Processing Planing plus Scheduling (right).}
    \label{fig:compar}
\end{figure}

\textbf{Novel Ingredient: }In our work, we model IPPS to equivalent MDP formulation, which is proved in Appendix A. While previous works on FJPS could only handle a single sequence, we designed a new heterogeneous graph based on the mentioned AND/OR graph to represent the state in IPPS, incorporating the Combination node category to extract information from the set of processing steps, and we use heterogeneous GNN to extract information and use the designed model to make decisions. To prevent the ineffective exploration costs caused by multiple production paths, we proposed an action space reduction based on the combination, ensuring that each trajectory selects exactly one processing sequence. Furthermore, we propose a reward function based on the estimated completion time for IPPS, and experiments show that it performs better than the naive reward function.
\begin{figure*}
    \centering
    \includegraphics[width=0.74\linewidth]{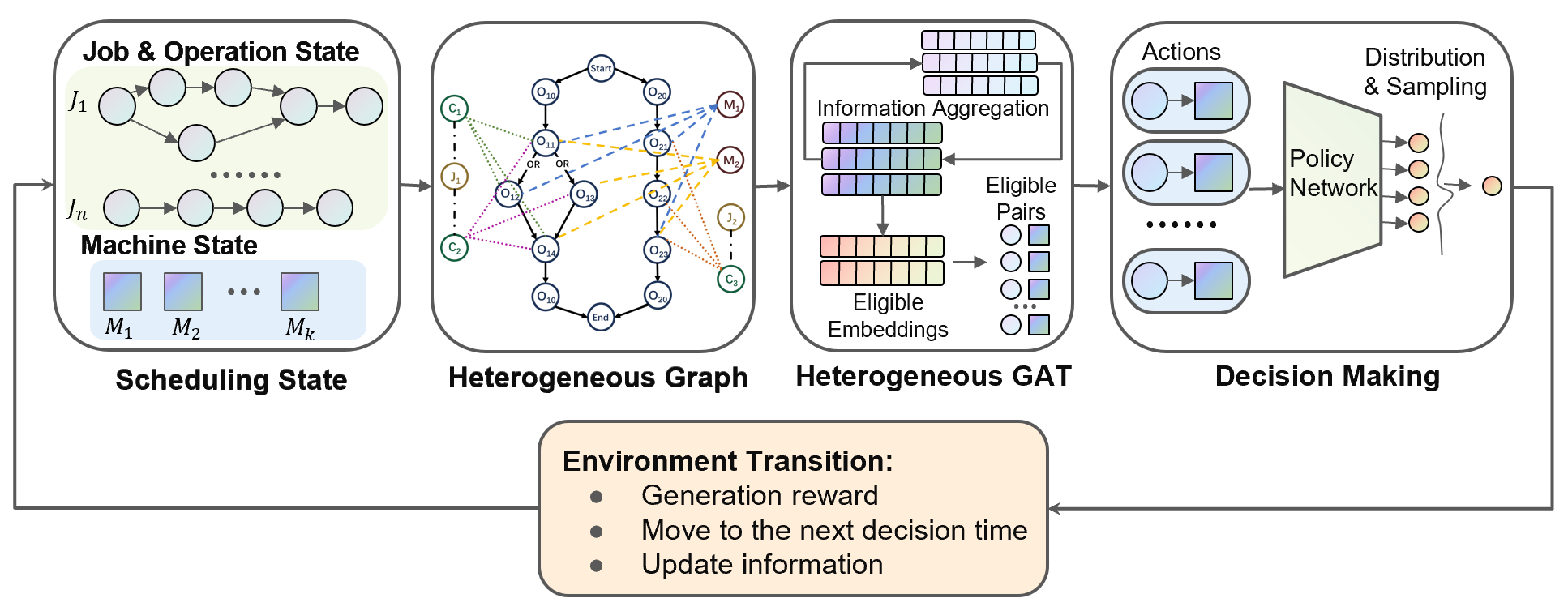}
    \caption{ The workflow of the proposed method.}
\end{figure*}

\section{MDP Formulation}
To define the Markov Decision Process (MDP) for the IPPS problem, we formulate the state space, action space, state transitions, and reward function as follows. At each decision point, following the completion of an operation, the agent extracts eligible operation-machine pairs based on the current state and decides which pair to perform next.

Notice that all decision points (or time-step) are at an operation's completion time. Whenever no operation can be processed or we choose to wait until the next time-step, we'll jump to the next time-step, which is the nearest completion time from the current time. The structure of the heterogeneous graph is illustrated as follows. We will prove the following theorem in Appendix A:

\textbf{Theorem 1} The introduced MDP formulation is equivalent to the original problem in terms of their optimal solutions.
\subsection{State}
At each step \( t \), the state \( s_t \) represents the status of IPPS problem job completion at the current time \( t \). We formulate the state using a dynamic heterogeneous graph, a directed structure composed of various types of nodes and edges, designed to extract information via GNNs. Previous studies on FJSP \cite{song2023,wang2023flexiblejobshopscheduling} usually include three types of nodes in heterogeneous graphs which cannot apply to IPPS problem directly since IPPS additionally requires selecting a combination following the choice of each OR-connector operation. Accordingly, we introduce nodes representing combinations and jobs within our graph to capture information from operations during GNN's information aggregation, which will be detailed in subsequent sections.

\begin{figure}[ht]
    \centering
    \includegraphics[width=0.3\linewidth]{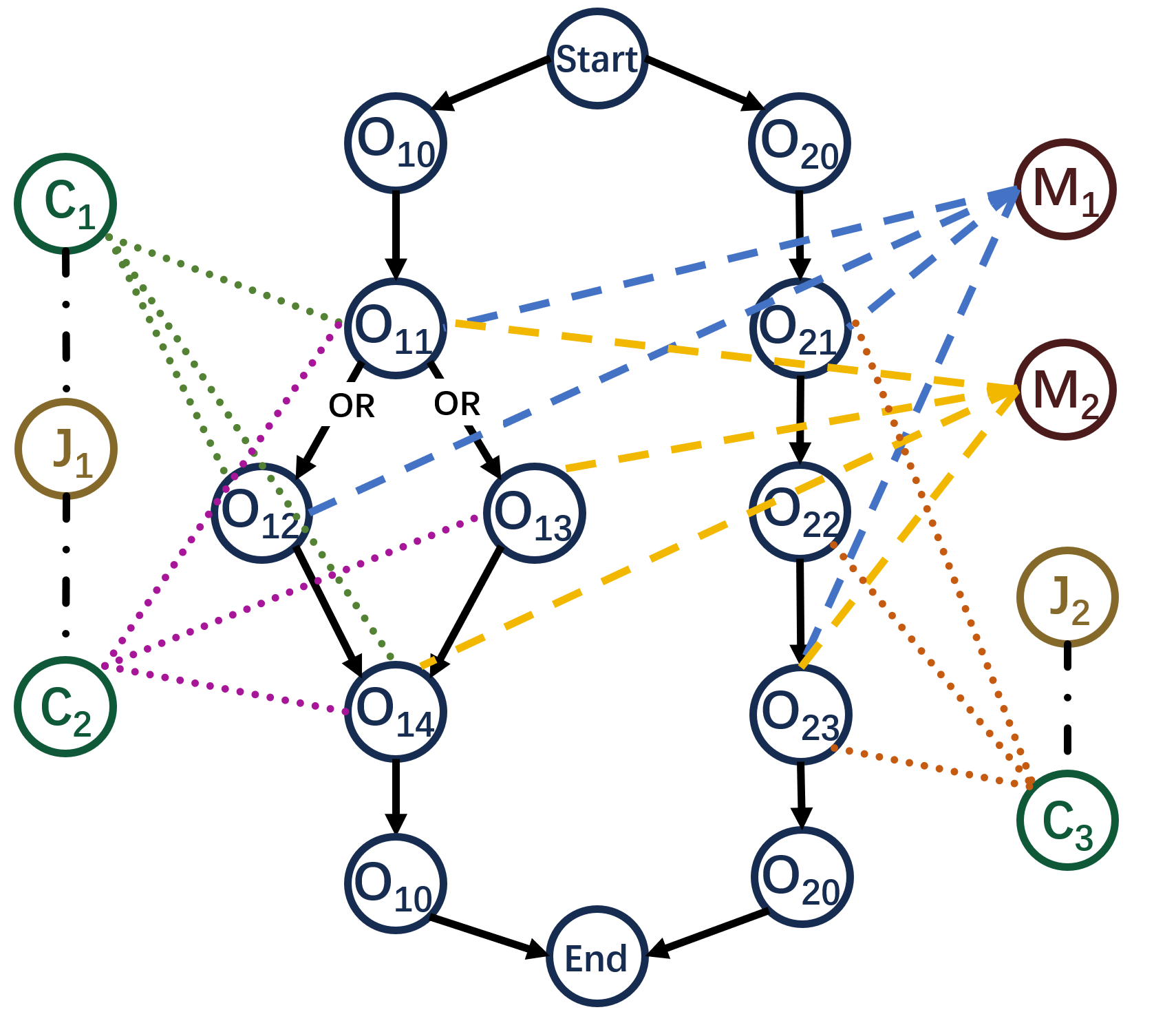}
    \caption{Example of heterogeneous graph.}
    \label{fig:simple-HGraph}
\end{figure}

\textbf{Node types:}
\begin{itemize}
    \item Operations (\(\mathcal{O}\)): \(O_{ij} \in \mathcal{O}\) denotes the \(i\)-th job's \(j\)-th operation. 
    \item Machines (\(\mathcal{M}\)): \(M_k \in \mathcal{M}\) denotes the \(k\)-th machines.
    \item Combinations (\(\mathcal{C}\)): \(C_{ih} \in \mathcal{C}\) denote a combination \(h\) of job \(i\), which is a set of operations. 
    \item Jobs (\(\mathcal{J}\)): \(J_i \in \mathcal{J}\) denotes the \(i\)-th job.
\end{itemize}

\textbf{Edge types:}
\begin{itemize}
    \item \(\mathcal{O} \longrightarrow\mathcal{O}\): represent the relationship of prerequisites, which means there is a directed edge from \(O_{ij}\) to \(O_{ij^{\prime}}\) if and only if \(O_{ij}\) is the last required predecessor operation of \(O_{ij^{\prime}}\).
    \item \(\mathcal{O} \longleftrightarrow\mathcal{C}\): there is an undirected edge between \(O_{ij}\) and \(C_{ih}\) when \(O_{ij} \in C_{ih}\).
    \item \(\mathcal{C} \longleftrightarrow\mathcal{J}\): there is an undirected edge between \(C_{ih}\) and \(J_{i}\) when \(C_{ih} \subset J_{i}\).
    \item \(\mathcal{O} \longleftrightarrow\mathcal{M}\): there is an undirected edge between \(O_{ij}\) and \(M_k\) when \(M_k\) can process \(O_{ij}\).
\end{itemize}
The details of the features in the  graph are shown in Appendix B. Figure~\ref{fig:simple-HGraph} is a heterogeneous graph of a simple IPPS problem (or initial state). Note that the "start" and "end" points are for representation purposes only, and the nodes marked \(O_{t0}\) are supernodes representing where the job \(J_i\) begins and ends with 0 processing time for all machines. 

\subsection{Action Space}

In this paper, an action \(a_t\) is defined as a feasible operation-machine pair (O-M pair), which means selecting an operation and a machine to process it. However, with numerous ineligible actions, there is a high cost to the agent's exploration. Therefore, we use two action space reduction strategies.

\textbf{Strategy I:}
The first strategy to reduce the action space is detecting feasibility, similar to previous FJSP works, such as \cite{song2023} and \cite{echeverria2024}. Specifically, candidates during the decision phase are restricted to O-M pairs that involve feasible operations and relevant idle machines, where feasible operations are those that have not been scheduled and all immediate prerequisite tasks have been completed. 

\begin{figure}[ht]
    \centering
    \includegraphics[width=0.5\linewidth]{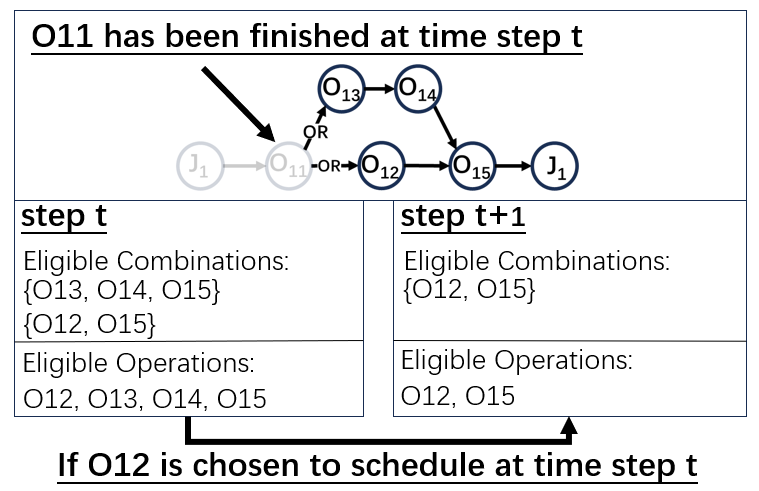}
    \caption{Example of action space reduction strategy II.}
    \label{fig:simple-OR}
\end{figure}

\textbf{Strategy II:}
\label{delete}
The IPPS problem has a more complex structure in which each job's production line is a directed acyclic graph. As described before, we conclude each job with several combinations, and it's evident that the number of combinations that might be chosen in the future, which we call "eligible", is non-increasing during dispatching. Inspired by this, our second reduction strategy is based on detecting eligible combinations. In more detail, we maintain \(C^*_t\), which is a set of eligible combinations. When one of the branches in the OR-connector is selected, all combinations containing other branches in this OR-connector will be removed from \(C^*_t\), and we will only consider operations that belong to at least one combination in \(C^*\). Figure\ref{fig:simple-OR} is a simple example that illustrates this strategy.

Above all, at time step \(t\), the action space \(\mathcal{A}_t\) is defined as:
\(\{waiting\} \cup \{(O_{ij}, M_k) | O_{ij} \in (\bigcup_{p} \{C^*_t\}_p) \cap R_t, M_k \in I_t \cap P_{ij}\}\)
, where \(\{C^*_t\}_p\) is the \(p\)-th component of \(C^*_t\) (set of eligible combinations), \(R_t\) is the set of feasible operations, \(I_t\) is the set of idle machines at \(t\), and \(P_{ij}\) is the set of machines that can process \(O_{ij}\).

\subsection{State Transitions}
In each time step, the heterogeneous graph is updated according to the following rules:

\begin{itemize}
\item Update the features in the graph as described previously.
\item Delete the finished operation nodes along with all arcs connected to these deleted operation nodes.
\item When a branch in an OR-connector is selected, delete the combination nodes and operation nodes according to the rule described in "action space reduction", as well as all arcs connected to these deleted nodes.
\item Delete machine nodes and job nodes that are no longer connected with any arcs.
\end{itemize}
\begin{figure}[ht]
    \centering
\includegraphics[width=0.5\linewidth]{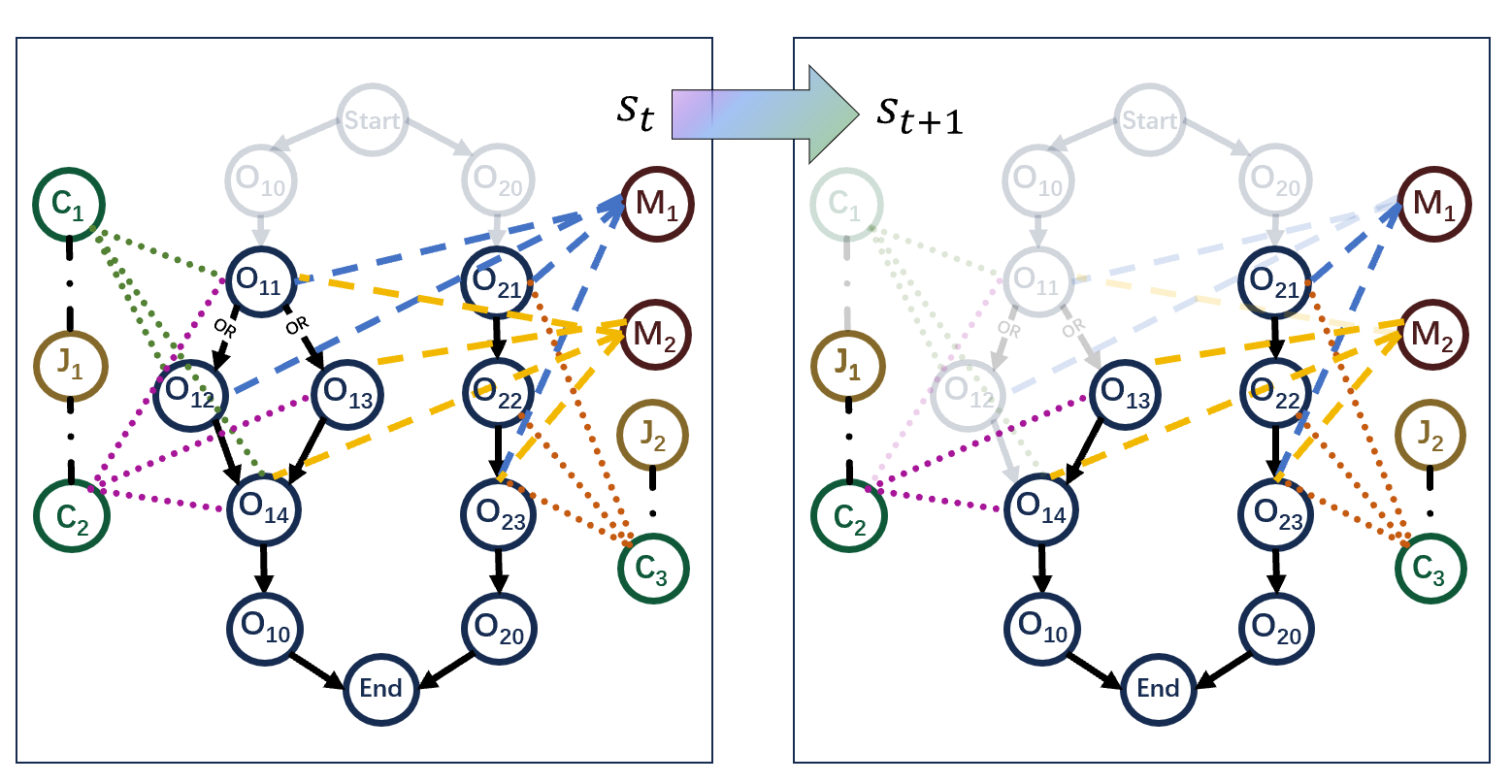}
    \caption{Example of state transition from time-step \(t\) to \(t+1\). }
    \label{fig:transition}
\end{figure}

Figure~\ref{fig:transition} illustrates a simple example of state transitions. At time step \(t\), \(O_{11}\) is scheduled but not yet completed, and \(M_1\) begins processing \(O_{12}\). Before \(t+1\), no O-M pair is selected because the only feasible operation, \(O_{21}\), can only be processed by \(M_1\), which is busy processing \(O_{11}\). Once \(O_{11}\) is completed, the corresponding node \(O_{11}\) and its related arcs are deleted. Post \(O_{11}\)'s completion, \(O_{12}\) and \(O_{13}\) become feasible, and \(M_1\) becomes idle. At \(t+1\), \(M_1\) is assigned to process \(O_{13}\), rendering \(C_1\) ineligible. Consequently, \(O_{12}\) and its related arcs are deleted. The features are updated accordingly. Note that super nodes do not impact the process as they finish immediately after being scheduled.

\subsection{Reward}
\label{reward}
The reward function is critical in the MDP, guiding the agent by providing feedback on the effectiveness of the current policy.

\textbf{Naive Reward:} Initially, we define a naive reward function \(r_0(s_t, a_t, s_{t+1}) = T(s_t) - T(s_{t+1})\), which represents the change in the actual completion time of the partial schedule. The cumulative naive reward without discount corresponds to the negative makespan. However, due to the significant variation among different IPPS problems, this type of cumulative reward fluctuates greatly, making it difficult to assess the current policy's performance.

\textbf{Suggested Reward:} To address this issue, we propose several advanced reward functions designed to evaluate policies based on relative improvement rather than absolute makespan. For initialization, we set the estimated processing time of each operation as \(\hat{p}_{ij} = \min_{k\in P_{ij}}{p_{ijk}}\) or \(\sum_{k\in P_{ij}}{p_{ijk}}/\|P_{ij}\|\), and aggregate these times by combinations to estimate their ending times \(\hat{T}_{ih}\). We then set \(\hat{T_i} = \min_{h \in J_{i}}{\hat{T}_{ih}}\) or \(\sum_{h \in J_i}{\hat{T}_{ih}}/\|J_i\|\) as the estimated ending time of jobs and \(\hat{T} = \max_i{\hat{T}_i}\) or \(\sum_{i \in J^*_t}{\hat{T}_i}/\|J^*_t\|\) as the estimated ending time of the entire problem, where \(J^*_t\) denotes jobs unfinished at \(t\). At each step, we update the estimated ending times of combinations and the problem using the real processing time and record the changes in estimated ending time as the reward \(r(s_t, a_t, s_{t+1}) = \hat{T}(s_t) - \hat{T}(s_{t+1})\). Our experiments demonstrate the effectiveness of this approach; see Appendix F.

\subsection{Policy}
The policy maps the state space to a probability distribution over the action space. For a given state \( s \), the policy \( \pi(s) \) provides the probabilities of the possible actions \( a \) in that state.

Mathematically, it can be expressed as:
\(
\pi : \mathcal{S} \rightarrow \mathcal{P}(\mathcal{A})
\)
where \( \mathcal{S} \) is the state space, \( \mathcal{P}(\mathcal{A}) \) is the probability distribution over the action space \( \mathcal{A} \).
We'll detail the policy design later.

\section{Neural Networks and Decision Making}
In this section,we introduce the method we use for feature extraction and algorithm for policy generation.
\subsection{Graph Representation}
DRL utilizes neural networks to extract state information. The states in our work are structured and have a variant number of nodes and edges, which cannot be represented as fixed-shaped tensors. Thus, we use graph neural networks(GNNs), which are shown to have higher performance in multiple fields with structured features, such as Combinatorial Optimization (CO) \cite{gasse2019exactcombinatorialoptimizationgraph,Liu_2024} and the recommendation system\cite{gao2023surveygraphneuralnetworks}.
One of the most popular and effective GNNs is graph attention networks (GATs), we apply GATv2\cite{gatv2} here, a modified version of GAT that shows higher efficiency in multiple fields. GATv2 is an aggregation-based GNN that uses dynamic attention mechanisms to aggregate information from neighbors connected by edges.
The attention score \(\alpha_{ij}\) is calculated as:
\begin{align*}
&\alpha_{i, j} = \frac{\exp(e_{ij})}{\sum_{k \in \mathcal{N}(i) \cup \{i\}} \exp(e_{ik})}, \text{ where } \\
&e_{ij} = 
\begin{cases} 
\mathbf{a}^{\top} \operatorname{LeakyReLU}(\boldsymbol{\Theta}_{s} \mathbf{x}_{i} + \boldsymbol{\Theta}_{t} \mathbf{x}_{j})\text{, if no edge features} \\ 
\mathbf{a}^{\top} \operatorname{LeakyReLU}(\boldsymbol{\Theta}_{s} \mathbf{x}_{i} + \boldsymbol{\Theta}_{t} \mathbf{x}_{j} + \boldsymbol{\Theta}_{e} \mathbf{e}_{ij})  \text{, else}
\end{cases}
\end{align*}

The new feature vector for node \(i\) is calculated as :
\begin{align*}
\mathbf{x}_{i}^{\prime}=\sum_{j \in \mathcal{N}(i) \cup\{i\}} \alpha_{i, j} \boldsymbol{\Theta}_{t} \mathbf{x}_{j}
\end{align*}

where \(\mathcal{N}(i)\) is a set of the first-hop neighbors of node \(i\), \(\mathbf{x}_i\) is input features(embedding) of node \(i\) and \(\mathbf{x}_i^\prime\) is output embedding,  \(\boldsymbol{\Theta}_s, \boldsymbol{\Theta}_t, \boldsymbol{\Theta}_e\) are independent linear transformations for features of start point, features of endpoint and edge features respectively. In our work, we use GATv2 in the heterogeneous graph by doing the aggregation steps described above for each edge type independently (the parameters are also independent) to extract embeddings, which we call HGAT.
\begin{figure}[ht]
    \centering
    \includegraphics[width=0.7\linewidth]{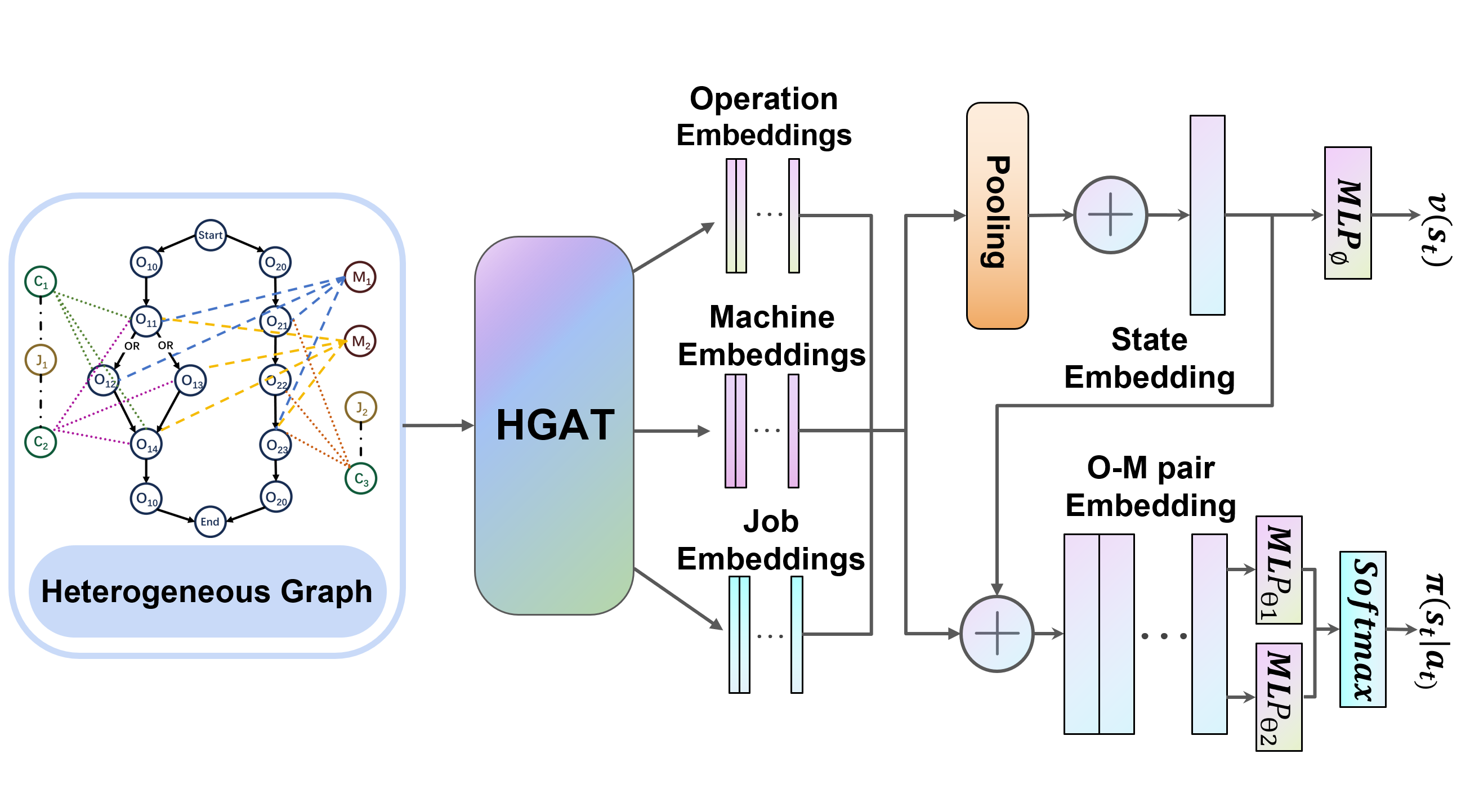}
    \caption{Architecture of neural network.}
    \label{fig:policy}
\end{figure}
\subsection{Policy Design}
A policy \(\pi_\theta(a_t|s_t)\)(or actor) is a mapping from state space \(\mathcal{S}\) to probability distribution over action space \(\mathcal{A}_t\).  
The Proximal Policy Optimization (PPO) algorithm \cite{schulman2017ppo} is one of our selected methods to train the agent, which also needs a value network \(v_\phi(s_t)\) (or critic). In our framework, the HGAT functions as an encoder, with its parameters being shared between the actor and critic albeit in distinct manners. The whole structure of our neural network is shown in Figure~\ref{fig:policy}.\par First, HGAT captures information from the heterogeneous graph and generates node embeddings for each type of nodes: \(\mathbf{h_{O_{ij}}}, \mathbf{h_{J_i}}, \mathbf{h_{M_k}}\). Note that node embeddings of combinations are not used in later steps. Define the set of feasible operations as \(R_t\) and \(F_t\) is the set of operations that are not scheduled and do not have unscheduled immediate prerequisite operations. \par To derive the action probabilities in the actor, let the set of eligible O-M pairs at time-step \(t\) be: \[ \mathbf{E} = \{(O_{ij},M_k)| O_{ij} \in (\bigcup_{p} \{C^*_t\}_p) \cap R_t , M_k \in I_t \cap P_{ij}\}\] To evaluate the waiting action, define future eligible O-M pairs as those with feasible operations and relevant machines processing other tasks: \[ \mathbf{E^\prime} = \{(O_{ij},M_k)| O_{ij} \in (\bigcup_{p} \{C^*_t\}_p) \cap F_t , M_k \in P_{ij}\}\]

The state embedding is: \[\mathbf{h_{t}}=\left[\frac{1}{|\mathcal{O}|} \sum_{i j } \mathbf{h_{O_{i j}}} \| \frac{1}{|\mathcal{M}|} \sum_{k} \mathbf{h_{M_{k}}}\| \frac{1}{|\mathcal{J}|} \sum_{i} \mathbf{h_{J_{i}}}\right]\] 

Also, we define the embedding of an O-M pair as:
\[
    \mathcal{H}_{ijk} = \left[ \mathbf{h_{O_{i j}}} \|\mathbf{h_{M_k}}\|\mathbf{h_{J_i}}\|\mathbf{h_t} \right]
\]
The priority \(\gamma\) of each eligible pair \((O_{ij}, M_k) \in \mathbf{E}\) is calculated as:\[
    \gamma(a_t|s_t)=\mathbf{MLP_{\theta_1}}(\mathcal{H}_{ijk})
\]
and the waiting priority \(\gamma\) is calculated as:\[
    \gamma(a_t|s_t) = \sum_{(O_{ij},M_k)\in \mathbf{E^\prime}}\alpha_{ijk}\mathbf{MLP_{\theta_2}}(\mathcal{H}_{ijk})
\]
where the attention score \(\alpha_{ijk}\) is calculated as follows:
\[
    \alpha_{ijk} = \frac{\operatorname{exp}\left[\mathbf{MLP_{\theta_1}}\left(\mathcal{H}_{ijk}\right)\right]}{\sum_{(O_{i^\prime j^\prime },M_k^\prime )\in \mathbf{E^\prime}}\operatorname{exp}\left[\mathbf{MLP_{\theta_1}}\left(\mathcal{H}_{i^\prime j^\prime k^\prime }\right)\right]}
\]
It should be noted that the independent parameters \(\theta_1\) and \(\theta_2\) are utilized to process eligible pairs and waiting actions, respectively. In addition, \(\theta_1\) is used to determine the weighting of each future eligible pair for the priority of waiting. This is because the importance of each future eligible pair varies, and \(\mathbf{MLP_{\theta_1}}\) can be used to assess the importance of each eligible pair, which in turn informs the weighting of future eligible pairs when calculating the priority of waiting actions.

The probabilities for each action \(a_t\) are calculated by a softmax layer over priorities, which is actually the actor \(\pi_\theta(a_t|s_t)\):
\[
    \pi_\theta(a_t|s_t) = \frac{\operatorname{exp}[\gamma(a_t|s_t)]}{\sum_{a^\prime_t\in \mathcal{A}_t}\operatorname{exp}[\gamma(a^\prime_t|s_t)]}
\]
The critic \(v_\phi(s_t)\) is derived from the state embedding using a MLP:\[
    v_\phi(s_t) = \mathbf{MLP}_\phi(\mathbf{h_t})
\]

\section{Experiments} 
\subsection{Dataset}
To evaluate the efficiency of our method, we conduct 3 experiments on two synthetic datasets of different scales and a public benchmark dataset from \cite{kim2003}.

\textbf{Synthetic Dataset: } Based on the structure of AND/OR graph, we generate jobs by computing direct acyclic graphs and add OR-paths recursively. After assigning machines and corresponding processing time to each operation under a certain distribution, we randomly combine jobs for various IPPS problems. For details, refer to Appendix D.

\textbf{Public Dataset: } A widely recognized benchmark in IPPS research introduced by \cite{kim2003}. Kim dataset contains 24 instances with different flexibility.
    
\subsection{Implementation Details}

Let \(N_1, N_2, N_3\) be the number of training, validation and testing instances, and \(m \times n\) denotes the problem size where \(m\) and \(n\) are the number of jobs and machines, respectively. Training instances are generated from 1000 jobs with corresponding machine numbers, which actually consist of \(N_1 = \binom{1000}{\text{job number}} \) problems.  We train three models, each based on synthetic datasets of different problem sizes: \(4 \times 5\), \(5 \times 3\), and \(6 \times 5\), which are all relatively small yet complicated and have a validation size of \(N_2 = 20\).

In Experiment 1, we evaluat the models on symmetric datasets that draw on the same size and distribution as the training set. To assess the generalization performance of the trained policies, Experiments 2 and 3 are conducted on a public dataset and on large-scale generated datasets with sizes \(16 \times 20\) and \(20 \times 25\), respectively. In all experiments, the size of all test sets is \(N_3 = 50\), and we use DRL-G (only choose actions with the highest probability) and DRL-S (choose actions according to the probability distribution) to evaluate policy performance.For details on infrastructure and configurations, please refer to Appendix E.

\subsection{Baseline}
We compare our model with two types of baseline algorithms.
\begin{itemize}
    \item \textbf{OR-Tools}: OR-Tools is an open source toolbox by Google \cite{ortools}, which includes the CP-SAT solver, a powerful constraint programming solver known for its strong performance in industrial scheduling \cite{dacoll2019industrial}. In this paper, we use the CP-SAT solver to implement the MILP model, which is shown in Appendix C, from \cite{jin2016} to obtain exact solutions in 1800 seconds (3600s for Experiment 3).
    \item \textbf{Greedy}: Dispatching rules have shown significant potential in IPPS problems \cite{ausaf2015pgha}. We select rules for operation sequencing and machine assignment, respectively. For operation selection, we use rules including Most Work Remaining (MWKR), Most Operations Remaining (MOR), First In First Out (FIFO), and probability described in \cite{ausaf2015pgha}. Machine assignment is based on the priority calculated from the Shortest Processing Time (SPT), Earliest Finish Time (EFT), and  Least Utilized Machine (LUM). In each problem, we randomly choose a combination for each job, reducing a complex IPPS problem to a simplified FJSP problem. Then, at each step, the next operation and the corresponding machine are determined by their rank under certain dispatching rules. Details are attached in Appendix G.
\end{itemize}

\subsection{Evaluation}
In order to compare method performance from both efficiency and quality, we use Makespan, Gap, and Time as evaluation metrics, where Makespan refers to the average final optimization objective, Gap is the average relative difference compared to Makespan of OR-Tools, and Time measures the average computation time. For the DRL-Sampling method, we run each instance 50 times and select the best value as the final result. For greedy algorithms, we choose the result of the best priority dispatching rules for each whole dataset.

\begin{table}[htbp]
\centering
\caption{Evaluation on instances of training sizes.}
\label{tab:exp1}
\medskip
\begin{tabularwithnotes}{cc c c ccc c}
{
\tnote[*]{The value in parentheses is the time for a single epoch. We ran multiple times and selected the Makespan.}
\tnote[\dag]{Gap is relative to OR-Tools.}
}
\toprule
\multicolumn{2}{c}{\textbf{Size}} & \textbf{OR-Tools} & \textbf{DRL-G} & \textbf{DRL-S} & \textbf{Best Greedy} \\
\midrule
\multirow{3}{*}{4$\times$5} & Makespan & 389.99 & 423.6 & \textbf{403.19} & 429.73 \\
 & Gap\textsuperscript{\dag}  & -      & 7.43\% & \textbf{2.92\%} & 8.79\% \\
 & Time (s) & 1769.31   & \textbf{7.65}   & 10.23(0.20)\tmark[*]    & 91.01(1.82)\tmark[*] \\
\midrule
\multirow{3}{*}{5$\times$3} & Makespan & 452.36 & 505.2 & \textbf{474.44} & 590.76 \\
 & Gap\textsuperscript{\dag}  & -      & 11.68\% & \textbf{4.88\%} & 30.60\% \\
 & Time (s) & 1789.79   & \textbf{7.90}   & 15.54(0.31)   & 48.45(0.97) \\
\midrule
\multirow{3}{*}{6$\times$5} & Makespan & 442.14 & 487.94 & \textbf{462.04} & 506.63 \\
 & Gap\textsuperscript{\dag}  & -      & 10.36\% & \textbf{4.50\%} & 14.59\% \\
 & Time (s) & 1798.84   & \textbf{8.37}   & 14.68(0.29)    & 139.76(2.80) \\
\bottomrule
\end{tabularwithnotes}
\end{table}

\begin{table}[htbp]
\centering
\renewcommand{\arraystretch}{0.8} 
\footnotesize
\caption{Performance on the Public Benchmark Kim dataset.}
\label{tab:kim_comparison}
\medskip
\begin{tabularwithnotes}{c c c c c}
{
\tnote[*]{LB indicates the lower bound on the makespan.}
\tnote[\dag]{Gap is relative to the OR-Tools makespan.}
}
\toprule
\textbf{Dataset} && \multicolumn{3}{c}{ \ \ \textbf{Kim}} \\
\cmidrule(lr){3-5}
& & \textbf{Makespan} & \textbf{Gap}\textsuperscript{\dag} & \textbf{Time(s)}  \\
\midrule
\textbf{LB}\textsuperscript{*} & & 387.93 & - & -  \\
\textbf{OR-Tools} & & 392.08 & 1.07\% &1796.37\\
\midrule
\textbf{DRL-G} & 4$\times$5 & 445.17 & 14.76\% &\textbf{11.94}\\
& 5$\times$3 & 451.33 & 16.34\% &17.79\\
& 6$\times$5 & 444.25 & 14.52\% &12.12\\
\textbf{DRL-S} & 4$\times$5 & 416.42 & 7.34\% &23.31(0.97)\\
& 5$\times$3 & 417 & 7.49\% &34.32(0.67)\\
& 6$\times$5 & \textbf{416.21} & \textbf{7.29}\% &23.87(0.48)\\
\textbf{Best Greedy} & & 458.38 & 18.16\% &82.24(1.64)\\
\bottomrule
\end{tabularwithnotes}
\end{table}

\subsection{Results and Analysis}
As shown in Table~\ref{tab:exp1}, in all training sizes, Experiment 1 demonstrates that our model consistently outperforms the Greedy method in terms of solution quality, while also being faster than all baseline methods. Our method consistently achieves a smaller gap in relation to OR-Tools compared to the Greedy method, with reductions ranging from 3\% to 25\%, highlighting the effectiveness of our approach.

Applying our model to the public dataset, the result of Experiment 2 in Table~\ref{tab:kim_comparison} shows similar advantages as Experiment 2, reflecting the generalization capability of our approach. Moreover, as shown by the detailed results of Kim dataset in Appendix H, our advantages get larger as the problem size increases.

Experiment 3 (see Table~\ref{tab:comparison}) further confirms our advantages in large-scale problems. Although OR-Tools cannot even find a feasible solution, the solving time of DRL is less than one minute. Also, DRL outperforms than OR-Tools in OR-Tools-feasible problems up to $11.35\%$. Strongly proving patterns learned in small and medium instances remain effective in solving larger problems. 

Due to page limitations, please refer to Appendix F for more numerical experiments.

\begin{table}[htbp]
\centering
\caption{Generalization performance on large-sized instances.}
\label{tab:comparison}
\begin{adjustbox}{width=0.7\textwidth}
\begin{tabularwithnotes}{lccccccc}
{
\tnote[*]{Problem set that is feasible to OR-Tools.}
\tnote[\ddag]{Gap is relative to OR-Tools.}
}
\toprule
\textbf{Dataset} & \textbf{Method} & & \multicolumn{2}{c}{\textbf{Makespan}} & \textbf{Gap}\textsuperscript{\ddag} & \textbf{Time(s)} \\
\cmidrule(lr){4-5}
& & & \textbf{Feasible}\textsuperscript{*} & \textbf{All}\textsuperscript{\dag} &  & \\
\midrule
\multirow{9}{*}{16$\times$20} & \textbf{OR-Tools} & & 566.07 & - & - & 1800 (7/50 time out) \\
\cmidrule(lr){2-7}
& \multirow{3}{*}{\textbf{DRL-G}} & 4$\times$5 & 538.37 & 537.84 & -4.89\% &\textbf{19.77}\\
& & 5$\times$3 &543.17 & 542.22 & -4.05\% & 32.67 \\
& & 6$\times$5 &542.2 & 539.54 & -4.22\% & 21.89 \\
\cmidrule(lr){2-7}
& \multirow{3}{*}{\textbf{DRL-S}} & 4$\times$5 & 518.22 & 527.06 & -8.45\% & \textbf{40.59(0.81)} \\
& & 5$\times$3 & \textbf{517.07} & \textbf{515.48} &\textbf{-8.66\%} & 65.15(1.30) \\
& & 6$\times$5 & 518.63 & 517.06 & -8.38\% & 41.00(0.82) \\
\cmidrule(lr){2-7}
& \textbf{Best Greedy} & & 544.80 & 544.89 & -3.76\% & 148.66(2.97) \\
\midrule
\multirow{9}{*}{20$\times$25} & \textbf{OR-Tools} & & 547.03 & - & - & 3600 (16/50 time out) \\
\cmidrule(lr){2-7}
& \multirow{3}{*}{\textbf{DRL-G}} & 4$\times$5 & 508.85 & 518.28 & -7.00\% & 26.87 \\
& & 5$\times$3 & 510.94 & 519.86 & -6.60\% & 42.37 \\
& & 6$\times$5 & 512.00 & 519.10 & -6.40\% & 28.23 \\
\cmidrule(lr){2-7}
& \multirow{3}{*}{\textbf{DRL-S}} & 4$\times$5 & \textbf{484.97} & \textbf{495.26} & \textbf{-11.35}\% & \textbf{56.34(1.13)} \\
& & 5$\times$3 & 486.15 & 496.84 & -11.13\% & 87.25(1.75) \\
& & 6$\times$5 & 485.67 & 495.24 & -11.21\% & 56.88(1.14) \\
\cmidrule(lr){2-7}
& \textbf{Best Greedy} & & 513.83 & 524.09 & -6.10\% & 186.42(3.73) \\
\bottomrule
\end{tabularwithnotes}
\end{adjustbox}
\end{table}

\section{Conclusions}
The IPPS problem is a critical challenge in combinatorial optimization in modern manufacturing systems. In this paper, we propose a novel DRL approach to solving IPPS, which results in a reduced exploration state space, faster solving speed, and smaller makespan gaps compared to greedy algorithms. Extensive experiments on synthetic and benchmark datasets demonstrate the efficiency of our method, particularly in large-scale problems.

\newpage

\bibliographystyle{unsrt}  
\bibliography{main}

\begin{thebibliography}{10}

\bibitem{SHAKYA2023120495}
Ashish~Kumar Shakya, Gopinatha Pillai, and Sohom Chakrabarty.
\newblock Reinforcement learning algorithms: A brief survey.
\newblock {\em Expert Systems with Applications}, 231:120495, 2023.

\bibitem{cappart2021combining}
Quentin Cappart, Thierry Moisan, Louis-Martin Rousseau, Isabeau Pr{\'e}mont-Schwarz, and Andre~A Cire.
\newblock Combining reinforcement learning and constraint programming for combinatorial optimization.
\newblock In {\em Proceedings of the AAAI Conference on Artificial Intelligence}, volume~35, pages 3677--3687, 2021.

\bibitem{zhang2020deep}
Rongkai Zhang, Anatolii Prokhorchuk, and Justin Dauwels.
\newblock Deep reinforcement learning for traveling salesman problem with time windows and rejections.
\newblock In {\em 2020 International Joint Conference on Neural Networks (IJCNN)}, pages 1--8. IEEE, 2020.

\bibitem{Zhang2023}
J.D. Zhang, Z.~He, W.H. Chan, and C.Y. Chow.
\newblock Deepmag: Deep reinforcement learning with multi-agent graphs for flexible job shop scheduling.
\newblock {\em Knowledge-Based Systems}, 259:110083, 2023.

\bibitem{Liu2023}
Y.~Liu, X.~Zuo, G.~Ai, and Y.~Liu.
\newblock A reinforcement learning-based approach for online bus scheduling.
\newblock {\em Knowledge-Based Systems}, 271:110584, 2023.

\bibitem{Liu2021}
Qihao Liu, Xinyu Li, Liang Gao, and Guangchen Wang.
\newblock Mathematical model and discrete artificial bee colony algorithm for distributed integrated process planning and scheduling.
\newblock {\em Journal of Manufacturing Systems}, 61:300--310, 2021.

\bibitem{Phanden2021}
Rakesh~Kumar Phanden, Ajai Jain, and J.~Paulo Davim, editors.
\newblock {\em Integration of Process Planning and Scheduling Approaches and Algorithms}.
\newblock Elsevier, 2021.

\bibitem{doi:10.1080/00207543.2021.1892230}
Gang~Du Yujie~Ma and Yingying Zhang.
\newblock Dynamic hierarchical collaborative optimisation for process planning and scheduling using crowdsourcing strategies.
\newblock {\em International Journal of Production Research}, 60(8):2404--2424, 2022.

\bibitem{samaranayake2012aircraft}
Premaratne Samaranayake and Sampath Kiridena.
\newblock Aircraft maintenance planning and scheduling: an integrated framework.
\newblock {\em Journal of Quality in Maintenance Engineering}, 18(4):432--453, 2012.

\bibitem{varela2021integrated}
Maria~LR Varela, Goran~D Putnik, Vijay~K Manupati, Gadhamsetty Rajyalakshmi, Justyna Trojanowska, and Jos{\'e} Machado.
\newblock Integrated process planning and scheduling in networked manufacturing systems for i4. 0: a review and framework proposal.
\newblock {\em Wireless Networks}, 27:1587--1599, 2021.

\bibitem{Zhou2019}
J.~Zhou, Y.~Zhou, B.~Wang, and J.~Zang.
\newblock Human-cyber-physical systems (hcpss) in the context of new-generation intelligent manufacturing.
\newblock {\em Engineering}, 5(4):624--636, 2019.

\bibitem{jin2016}
Liangliang Jin, Qiuhua Tang, Chaoyong Zhang, Xinyu Shao, and Guangdong Tian.
\newblock More milp models for integrated process planning and scheduling.
\newblock {\em International Journal of Production Research}, 54:4387 -- 4402, 2016.

\bibitem{barzanji2020}
R.~Barzanji, B.~Naderi, and M.~A. Begen.
\newblock Decomposition algorithms for the integrated process planning and scheduling problem.
\newblock {\em Omega}, 93, 2020.

\bibitem{kim2003}
Y.~K. Kim, K.~Park, and J.~Ko.
\newblock A symbiotic evolutionary algorithm for the integration of process planning and job shop scheduling.
\newblock {\em Computers \& Operations Research}, 30(8):1151--1171, 2003.

\bibitem{li2018}
Y.~Li, L.~Gao, Q.~K. Pan, L.~Wan, and K.-M. Chao.
\newblock An effective hybrid genetic algorithm and variable neighborhood search for integrated process planning and scheduling in a packaging machine workshop.
\newblock {\em IEEE Transactions on Systems, Man, and Cybernetics: Systems}, 49(10):1933--1945, 2018.

\bibitem{luo2021}
S.~Luo, L.~Zhang, and Y.~Fan.
\newblock Dynamic multi-objective scheduling for flexible job shop by deep reinforcement learning.
\newblock {\em Computers \& Industrial Engineering}, 159:107489, 2021.

\bibitem{song2023}
Wen Song, Xinyang Chen, Qiqiang Li, and Zhiguang Cao.
\newblock Flexible job-shop scheduling via graph neural network and deep reinforcement learning.
\newblock {\em IEEE Transactions on Industrial Informatics}, 19(2):1600--1610, 2023.

\bibitem{echeverria2023}
I.~Echeverria, M.~Murua, and R.~Santana.
\newblock Solving the flexible job-shop scheduling problem through an enhanced deep reinforcement learning approach.
\newblock {\em arXiv preprint arXiv:}, 2023.

\bibitem{li2007}
W.~D. Li and C.~A. McMahon.
\newblock A simulated annealing-based optimization approach for integrated process planning and scheduling.
\newblock {\em International Journal of Computer Integrated Manufacturing}, 20(1):80--95, 2007.

\bibitem{leung2010}
C.~Leung, T.~Wong, K.-L. Mak, and R.~Y. Fung.
\newblock Integrated process planning and scheduling by an agent-based ant colony optimization.
\newblock {\em Computers \& Industrial Engineering}, 59(1):166--180, 2010.

\bibitem{li2010}
X.~Li, X.~Shao, L.~Gao, and W.~Qian.
\newblock An effective hybrid algorithm for integrated process planning and scheduling.
\newblock {\em International Journal of Production Economics}, 126(2):289--298, 2010.

\bibitem{wu2021}
X.~Wu and J.~Li.
\newblock Two layered approaches integrating harmony search with genetic algorithm for the integrated process planning and scheduling problem.
\newblock {\em Computers \& Industrial Engineering}, 2021.

\bibitem{han2021}
B.~Han and J.~Yang.
\newblock A deep reinforcement learning based solution for flexible job shop scheduling problem.
\newblock {\em International Journal of Simulation Modelling}, 20(2):375--386, 2021.

\bibitem{zhang2020}
C.~Zhang, W.~Song, Z.~Cao, J.~Zhang, P.~S. Tan, and X.~Chi.
\newblock Learning to dispatch for job shop scheduling via deep reinforcement learning.
\newblock {\em Advances in Neural Information Processing Systems}, 33, 2020.

\bibitem{ho1996solving}
YC~Ho and CL~Moodie.
\newblock Solving cell formation problems in a manufacturing environment with flexible processing and routing capabilities.
\newblock {\em International Journal of Production Research}, 34(11):2901--2923, 1996.

\bibitem{wang2023flexiblejobshopscheduling}
Runqing Wang, Gang Wang, Jian Sun, Fang Deng, and Jie Chen.
\newblock Flexible job shop scheduling via dual attention network based reinforcement learning, 2023.

\bibitem{echeverria2024}
I.~Echeverria, M.~Murua, and R.~Santana.
\newblock Leveraging constraint programming in a deep learning approach for dynamically solving the flexible job-shop scheduling problem.
\newblock {\em arXiv preprint arXiv:2403.09249}, 2024.

\bibitem{gasse2019exactcombinatorialoptimizationgraph}
Maxime Gasse, Didier Chételat, Nicola Ferroni, Laurent Charlin, and Andrea Lodi.
\newblock Exact combinatorial optimization with graph convolutional neural networks, 2019.

\bibitem{Liu_2024}
Tianhao Liu, Shanwen Pu, Dongdong Ge, and Yinyu Ye.
\newblock Learning to pivot as a smart expert.
\newblock {\em Proceedings of the AAAI Conference on Artificial Intelligence}, 38(8):8073–8081, March 2024.

\bibitem{gao2023surveygraphneuralnetworks}
Chen Gao, Yu~Zheng, Nian Li, Yinfeng Li, Yingrong Qin, Jinghua Piao, Yuhan Quan, Jianxin Chang, Depeng Jin, Xiangnan He, and Yong Li.
\newblock A survey of graph neural networks for recommender systems: Challenges, methods, and directions, 2023.

\bibitem{gatv2}
Shaked Brody, Uri Alon, and Eran Yahav.
\newblock How attentive are graph attention networks?
\newblock {\em ArXiv}, abs/2105.14491, 2021.

\bibitem{schulman2017ppo}
John Schulman, Filip Wolski, Prafulla Dhariwal, Alec Radford, and Oleg Klimov.
\newblock Proximal policy optimization algorithms.
\newblock {\em arXiv preprint arXiv:1707.06347}, 2017.

\bibitem{ortools}
{Google}.
\newblock Or-tools.
\newblock \url{https://github.com/google/or-tools}, 2024.
\newblock Accessed: 2024-07-30.

\bibitem{dacoll2019industrial}
G.~Da Col and E.~C. Teppan.
\newblock Industrial size job shop scheduling tackled by present day cp solvers.
\newblock In {\em International Conference on Principles and Practice of Constraint Programming}, pages 144--160. Springer, 2019.

\bibitem{ausaf2015pgha}
Muhammad~Farhan Ausaf, Liang Gao, Xinyu Li, and Ghiath~Al Aqel.
\newblock A priority-based heuristic algorithm (pbha) for optimizing integrated process planning and scheduling problem.
\newblock {\em Cogent Engineering}, 2(1):1070494, 2015.

\bibitem{pytorch}
Adam Paszke, Sam Gross, Francisco Massa, Adam Lerer, James Bradbury, Gregory Chanan, Trevor Killeen, Zeming Lin, Natalia Gimelshein, Luca Antiga, et~al.
\newblock Pytorch: An imperative style, high-performance deep learning library.
\newblock {\em Advances in neural information processing systems}, 32, 2019.

\bibitem{pyg}
Matthias Fey and Jan~Eric Lenssen.
\newblock Fast graph representation learning with pytorch geometric.
\newblock {\em arXiv preprint arXiv:1903.02428}, 2019.

\end{thebibliography}
\begin{appendices}

\section{Properties of MDP Formulation}

Although relevant DRL algorithms proposed before have reached competitive performance, they all neglect the value of waiting
and arrange feasible ope-ma pairs as soon as possible, which makes them impossible to reach optimality in some cases. Initially, we present a straightforward example to demonstrate this and subsequently show that our formulation is equivalent to the original IPPS problem.
\subsection{Example}
There's a simple example (FJSP or IPPS) explaining why the policy without waiting may miss the (near) optimal solution.\\
\begin{table}[H]
  \centering
  \begin{tabular}{c|cc}
      \toprule
      & ma1 & ma2 \\
      \midrule
      ope1 & 1 & 1 \\
      ope2 & 3 & 1 \\
      ope3 & 4 & 2 \\
      \bottomrule
  \end{tabular}
\end{table}
In this problem, we have only 2 jobs, 2 machines and 3 operations. Operations for job 1 and 2 is \{ope1, ope2\} and \{ope3\} respectively and their
processing times are above. The optimal schedule for this problem is \{(ope1,ma1),(ope3,ma2),(ope2,ma2)\}, which has a waiting decision and makespan of 3.
However, in the policy without waiting, we would arrange an ope-ma pairly if there exists feasible pairs.
Hence, the algorithm would process ope1 and ope3 at $t=0$ and job1 would be available again at $t=1$. Due to the loss of waiting and ope3's long processing times,
 ope2 would be arranged to the same machine as ope1, which means that each job can only be processed on a single machine. Because the total processing times of two jobs on ma1 are both 4, the best makespan under algorithm without waiting would be 4 and cannot reach the optimality.
\begin{figure}[ht]
    \centering
    \begin{subfigure}{.25\textwidth}
        \centering
        \includegraphics[width=\linewidth]{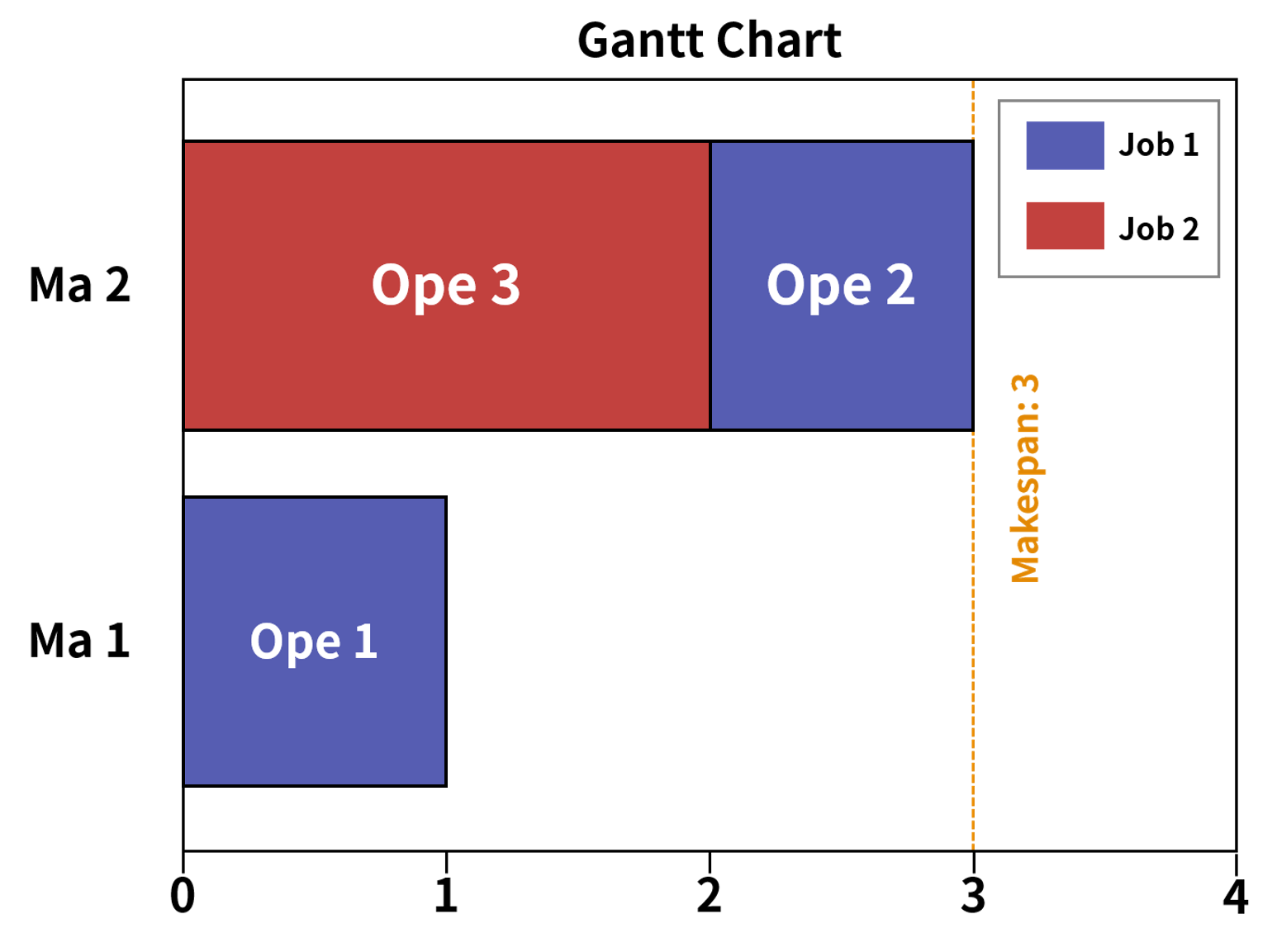}
        \caption{Optimal}
    \end{subfigure}%
    \begin{subfigure}{.25\textwidth}
        \centering
        \includegraphics[width=\linewidth]{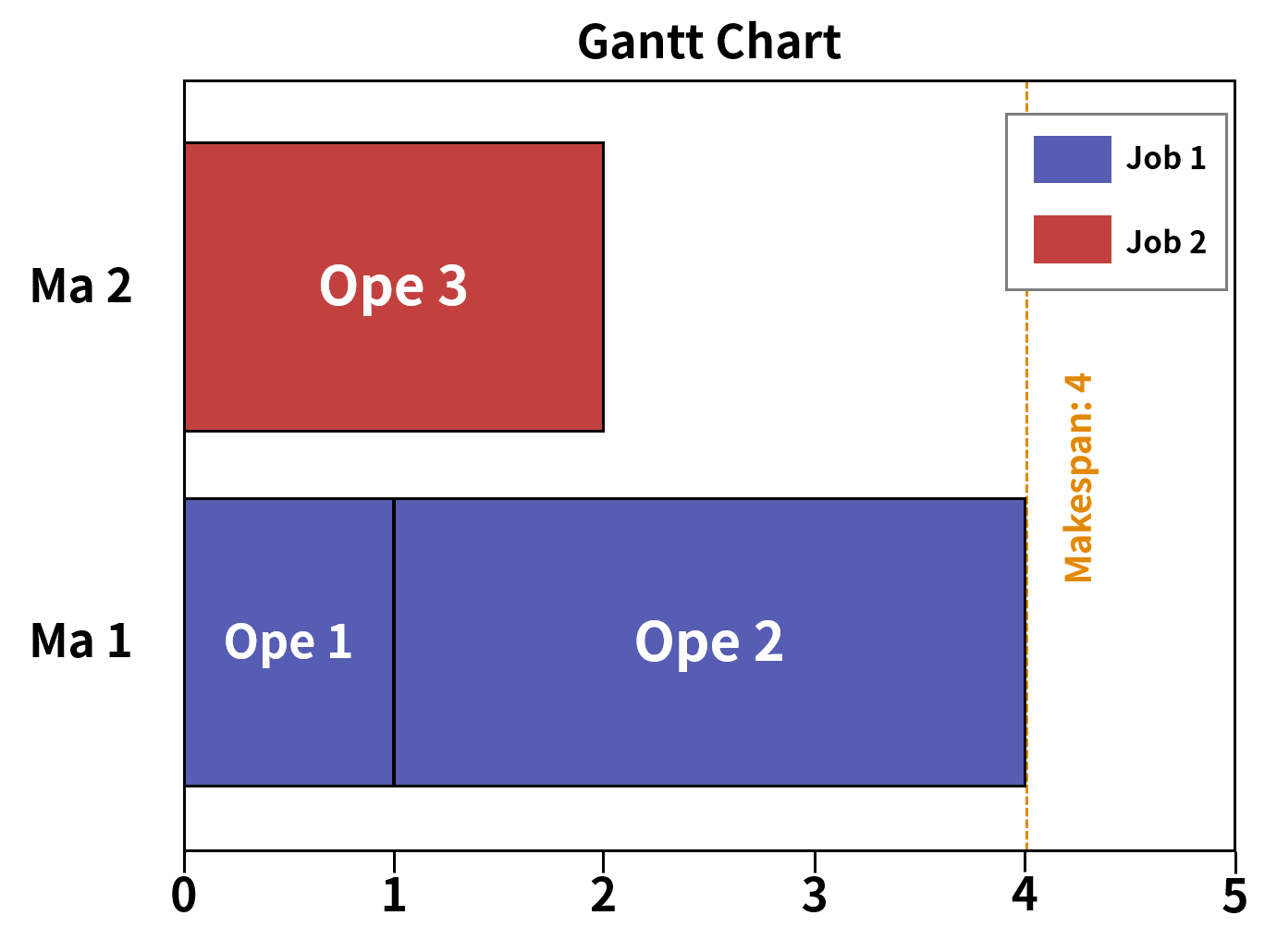}
        \caption{No waiting}
    \end{subfigure}
    \caption{}
\end{figure}
\subsection{Proof of Theorem 1 (Equivalence)} 
Here we prove that our MDP is equivalent to the original IPPS problem.
We define an IPPS solution as \(D = \{d, s, e\}\), where \(d = \{d_i | i=1,2,\dots,|d|\}\) is the set of schedules, with \(d_i\) being the \(i\)-th schedule allocating operation \(O_i\) belonging to \(J_i\) to machine \(M_i\). The functions \(s(d_i)\) and \(e(d_i)\) represent the start and end times of \(d_i\), and the makespan \(M=\max\{e(d_i)|d_i\in d\}\). We can assume \(s(d_i)\)  are non-decreasing with \(i\). 

\textbf{Recall}: A schedule \(d_i = (O_i, M_i, J_i)\) is eligible \(\iff\) all prerequisite operations of operation \(O_i\) have been finished, all operations that belong to job \(J_i\) and are not prerequisite for \(O_i\) hasn't been scheduled and machine \(M_i\) and job \(J_i\) are idle, which is equivalent to:
\begin{align*}
    &\forall O_m\in \mathcal{P}(O_i)\   \exists d_j \text{ s.t. } O_j \in d_i \And e(d_j) \leq s(d_i)\\
    &\forall O_m \in \mathcal{A}(O_i)\   \forall d_j\in\{d_p|s(d_p) \leq s(d_i)\},   O_m \not= d_j\\
    & e(d_t) \leq s(d_i), \forall d_t \text{ s.t } M_i = M_t,t\leq i\\
    & e(d_t) \leq s(d_i), \forall d_t \text{ s.t } J_i = J_t,t\leq i
\end{align*}
where \(\mathcal{P}(O_i)\) is the set of prerequisite operations of operation \(O_i\), and \(\mathcal{A}(O_i)\) is the set of operations belong to same job with \(O_i\) and not in \(\mathcal{P}(O_i)\). For convenient discussion, we define \(\mathcal{R}_i=\{(O_h, M_h, J_h) | O_h=O_i \bigvee M_h=M_i\}\) for \(d_i = (O_i, M_i, J_i) \) as the set of schedules containing   prerequisite operations of \(O_i\), \(M_i\) and \(J_i\) and the condition above will be:
\begin{align*}
    &\forall O_m\in \mathcal{P}(O_h)\   \exists d_j\in \mathcal{R}_i \text{ s.t. } O_m = O_j \tag{1}\\
    & e(d_t) \leq s(d_i), \forall d_t \in \mathcal{R}_i \tag{2}\\
    &\forall O_m \in \mathcal{A}(O_h)\   \forall d_j\in\{d_p|s(d_p) \leq s(d_i)\},   O_m \not= O_j \tag{3}
\end{align*}

\textbf{Lemma 1:} for \(\forall\) IPPS solution \(D\) with makespan \(M\), \(\exists D^\prime\) with makespan\(M^\prime\) that satisfies:
\begin{align*}
    &\forall i\ \  \exists j \text{ s.t. } s(d_i) = e(d_j) \tag{4}\\ &M^\prime \leq M\tag{5}
\end{align*}

\textbf{Proof:} Consider an IPPS solution \(D\) that does not satisfy (4). Assume \(d_i\) is the first schedule that $\exists! j \ \text{s.t. } s(d_i) = e(d_j) $. Notice that all schedules in a solution are eligible, we have:\[
    \forall d_j \in \mathcal{R}_i, e(d_j)<s(d_i)
\]
We create \(D^\prime\) from \(D\) by following rule:
\begin{itemize}
    \item \(d^\prime = d\) 
    \item \(s^\prime(d^\prime_i) = \min \{e(d_j)|d_j \in \mathcal{R}_i\} < s(d_i)\)
    \item \(e^\prime(d^\prime_i) = s(d^\prime_i) + \text{processing time of } d_i < e(d_i)\)
    \item \(s^\prime(d^\prime_t) = s(d_t), e^\prime(d^\prime_t) = e(d_t), \forall t \not= i\)
\end{itemize}
Then, \(D^\prime\) is also a solution because all schedules also satisfy (1) and (2) and \(d^\prime_i\) satisfies (4). (1) holds because \(d^\prime = d\) . And the reason that (2)(3) holds for all schedules in  \(D^\prime\) is:
\begin{itemize}
    \item From our construction, \(d^\prime_i\) satisfies (2).
    \item \(d^\prime_i\) satisfies (3), because \[\{d_p|s^\prime(d_p) \leq s^\prime(d_i)\}\subseteq\{d_p|s(d_p) \leq s(d_i)\}\]
    \item For \(t \not= i\), \(d^\prime_t\) satisfies (2) because:
    \begin{itemize}
        \item if \( \forall j \ d^\prime_j \in \mathcal{R}_t\), \(e^\prime(d^\prime_j) \leq e(d_j) \leq s(d_t) = s^\prime(d^\prime_t)\) 
        \item else, \(e^\prime(d^\prime_j) = e(d_j) \leq s(d_t) = s^\prime(d^\prime_t)\) 
    \end{itemize} 
    \item For \(t \not= i\), \(d^\prime_t\) satisfies (3) because (assume the operation in \(d^\prime_i, d^\prime_t\) are \(O_m, O_n\) respectively):
    \begin{itemize}
        \item if \(O_m\not\in \mathcal{A}(O_n)\), (3) holds as everything doesn't change.
        \item else, we only need consider \(d_i\) as:
        \[d^\prime_j = d_j, s^\prime(d^\prime_j) = s(d_j), e^\prime(d^\prime_j) = e(d_j), j\not=i\]
        we know \(s(d_t)<s(d_i)\) from (3) and \(d_t\in \mathcal{R}_i\), which means 
        \[ s^\prime(d^\prime_t) < e^\prime(d^\prime_t) \leq s^\prime(d^\prime_i)\]
        So, (3) holds
    \end{itemize}
\end{itemize}
Also, as \(e^\prime(d^\prime_j) \leq e(d_j), \forall j\), \(M^\prime \leq M\).
Keep following the procedure described above, we can find a solution \(\hat D\) satisfies (4)(5).

\textbf{Lemma 2:} For any solution \(D\) satisfying (4), it can be attained by MDP.

\textbf{Proof:} 
In our MDP formulation described before, a schedule \(d_i = (O_i, M_i, J_i)\) is eligible in the time step \(t\) if (1) (2) (3) holds when \(s(d_i) = t\), so we need to prove that all \(s(d_i)\) can be attained in MDP. Notice that the support set of \(\{s(d_j)|d_j\in D\}\) can be written as a non-decreasing sequence \(0 = s_0, s_1, \dots s_N\). We prove by induction:\\
(1)Basic case, start time 0 can be attained by MDP and MDP can schedule as same as \(D\) for \(d_i \text{ with } s(d_i)=0\), which is obvious.\\
(2) Assume \(\forall s_i \in \{s_t|s_t\leq s_k, k<N\}\) can be attained by MDP, and MDP have same schedules with start time \(\leq s_k\) with original solution \(D\).Notice We must has a completion time that equals \(s_{k+1}\) as all schedules with start time \(\leq s_k\) is as same as \(D\) and we can keep waiting until time equals \(s_{k+1}\) because each time-step is a completion time. Then from (4), we can know \(S_{k+1}\) can be attained and MDP have same schedules with start time \(\leq s_k+1\) with original solution \(D\).

A trajectory given from MDP must be a feasible solution because the action space design guarantee (1)(2)(3) must hold. Also, from Lemma 1 and Lemma 2, we know that for any optimal solution in original problem, we can get another optimal solution that can be given by MDP by a simple adjustment. So, it's equivalent to optimize original problem and MDP formulation problem.

\section{Graph Features}
Raw features for Nodes:
\begin{itemize}
    \item Operations:
        \begin{itemize}
            \item Number of prerequisite operations.
            \item Whether this operation is scheduled.
            \item Whether operation is feasible.
            \item Waiting time for a feasible operation(=0 for infeasible operations).
            \item Time remained for the processing operation that is scheduled(=0 for unscheduled operation).
        \end{itemize}
    \item Machines: \begin{itemize}
            \item Number of neighboring operations.
            \item The time when this machine becomes idle.
            \item Utilization: ratio of the non-idle time to the total production time.
            \item Whether the machine is working.
            \item Idle time from available time to now.
            \item Time remained for working machine to finish the working operation.
        \end{itemize}
    \item Combinations: \begin{itemize}
            \item Estimated completion time for these combinations.
            \item Ratio of estimated end time to estimated min end time for combinations associated to the same job.
        \end{itemize}
    \item Jobs: \begin{itemize}
            \item Ratio of estimated end time to max estimated end time for all jobs.
        \end{itemize}

\end{itemize}
Raw features for edges:
\begin{itemize}
    \item Operations \(\Leftrightarrow\) Machines:
        \begin{itemize}
            \item Processing time used by the machine to process the operation. 
        \end{itemize}
\end{itemize}
\label{feature}

\section{MILP Model}
To attain the optimal solution, we implement the MILP (Mixed Integer Linear Programming) 
model suggested by\cite{jin2016} using the CP-SAT solver of OR-Tools, an open source optimization library developed by Google. For more details, please refer cited paper.\\
\subsection*{Signal:}
\begin{itemize}
    \item \( i, i' \): Job indices, \(1 \leq i \leq |n|\)
    \item \( j, j' \): Operation indices, \(1 \leq j \leq |n_i|\)
    \item \( k, k' \): Machine indices
    \item \( h \): Combination index
    \item \( O_{ij} \): The \(j\)th operation of the \(i\)th job
    \item \( O_{ihj} \): The \(j\)th operation of the \(h\)th combination of the \(i\)th job
\end{itemize}

\subsection*{Sets and Parameters:}
\begin{itemize}
    \item \( A \): A very large positive number.
    \item \( P_{ijk} \): Processing time of job \(O_{ij}\) on machine \(k\).
    \item \( R_h \): The set of jobs contained in the \(h\)th combination of job \(i\).
    \item \( K_i \): The set of all combinations of job \(i\).
    \item \( n \): The set of jobs.
    \item \( n_i \): The set of all jobs in the network diagram of job \(i\).
    \item \( M_{ij} \): The set of all candidate machines for job \(O_{ij}\).
    \item \( Q_{ijj'} \): Equals 1 if job \(O_{ij}\) is processed  ly or directly before job \(O_{ij'}\); equals 0 otherwise.
\end{itemize}

\subsection*{Variables:}

\begin{itemize}
    \item \(X_{ihjk}\): Equals 1, if Operation \(O_{ihj}\) is allocated to Machine \(k\). 
    \item \(C_{ihj}\): The completion time of the operation \(O_{ihj}\).
    \item \( C_{\max} \): Maximum makespan.
    \item \( Y_{ih} \): Equals 1 if the \(h\)th combination of job \(i\) is selected; equals 0 otherwise.
    \item \( Z_{ijj'} \): Equals 1 if job \(O_{ij}\) is processed directly or indirectly
    before job \(O_{ij'}\); equals 0 otherwise.
\end{itemize}

\subsection*{Constraints:}
For each job \( i \), exactly one combination \( h \) must be selected.
\begin{align*}
    \sum_{h \in K_{i}} Y_{ih} = 1, \quad \forall i
\end{align*}

If combination \( h \) for job \( i \) is selected, the sum of operations in the combination is equal to 1, else it equals 0.
\begin{align*}
    \sum_{k \in \mathcal{M}_{ij}} X_{ihjk} = Y_{ih}, \quad \forall i, h \in K_{i}, j \in R_{ih}
\end{align*}

If combination \( h \) is not selected, the completion time of all operations in the combination is 0.
\begin{align*}
    A \cdot Y_{ih} \geq C_{ihj}, \quad \forall i, h \in K_i, j \in R_{ih}
\end{align*}

Guarantee correct completion time of operations.
\begin{align*}
    C_{ihj'} \geq C_{ihj} + \sum_{k' \in M_{y'}} X_{ihj'k'} P_{ij'k'}, \\\quad \forall i, h \in K_{i}, \; j, j' \in R_{ih}, \; j \neq j', \; V_{ij'} = 1
\end{align*}

\begin{align*}
    Z_{ijj'} + Z_{ij'j} = 1, \quad \forall i, j, j' \in n_{i}, \; Q_{ijj'} + Q_{ij'j} = 0, \; j \neq j'
\end{align*}

\begin{align*}
    C_{ihj'} \geq C_{ihj} + \sum_{k' \in M_{v}} X_{ihj'k'} P_{ij'k'} - A(1 - Z_{ijj'}), \\\quad \forall i, h \in K_{i}, \; j, j' \in R_{ih}, \; j \neq j'
\end{align*}

Guarantee the processing order of operations in the same machines.
\begin{align*}
    &C_{ih'j'} \geq C_{ihj} + X_{i'h'j'k'} P_{i'k'}\\&\quad -  A(1 - W_{ijj'j'}) - A(2 - X_{ihjk} - X_{i'h'j'k'}),\\ 
    &C_{ihj} \geq C_{i  h' j'} + X_{ihjk} P_{ijk} - A \cdot W_{ijj'j'}\\ &\quad  - A(2 - X_{ihjk} - X_{ih'j'k'}), \\
    &\quad \forall i, i', \; i \neq i', \; h \in K_{i}, \; h' \in K_{i'}, \;j \in R_{ih}, \; \\&\quad j' \in R_{i'h'}, \; k, k' \in M_{ij} \cap M_{ij'}, \; k = k'
\end{align*}

Makespan is the maximum completion time of all operations.
\begin{align*}
    C_{\max} \geq C_{ihj}, \quad \forall i, h \in K_{i}, \; j \in R_{ih}
\end{align*}

\section{Problem Generation Approach}
 Based on the structure of AND/OR graph, we generate jobs by computing direct acyclic graphs and add OR-paths recursively. After assigning machines and corresponding processing time to each operation under a certain distribution, we randomly combine jobs for various IPPS problems. 
 
Here, we define OR-link path as the sub-path, a path that is not the OR-link path as the main-path.
\subsection{Non-Conforming Situations}
\begin{figure}[ht]
    \centering
    \begin{subfigure}{0.4\linewidth}
        \centering
        \includegraphics[width=0.9\linewidth]{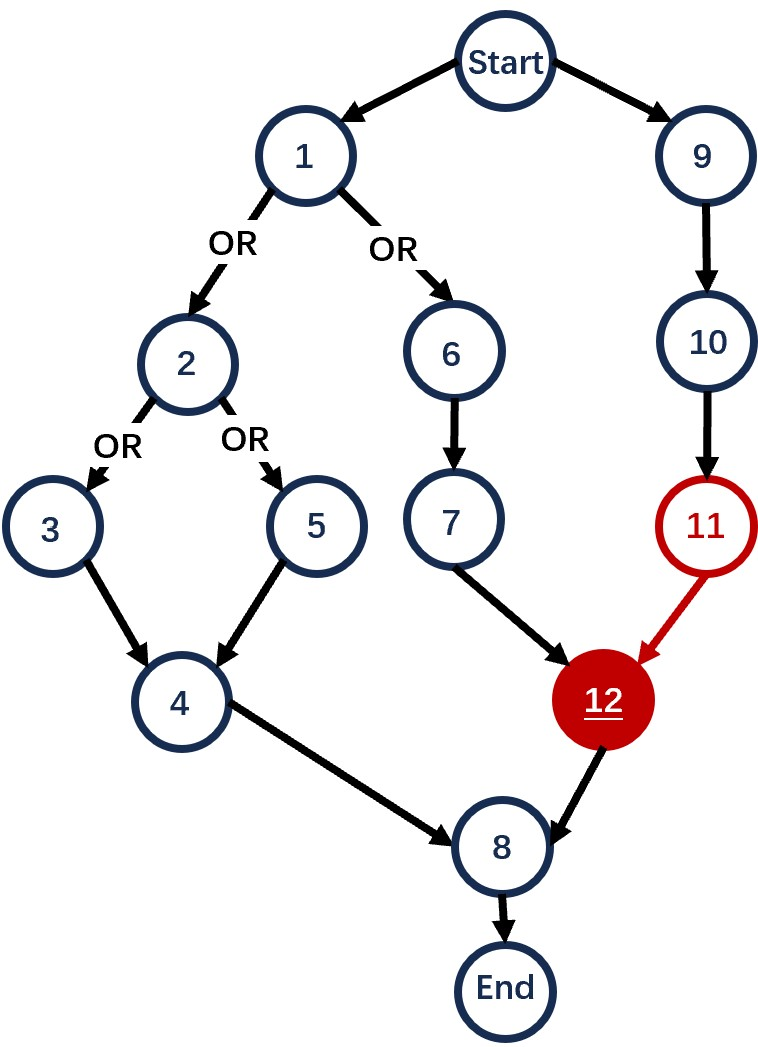}
        \caption{Main path to sub-path}
        \label{fig:main_to_sub}
    \end{subfigure}
    \begin{subfigure}{0.41\linewidth}
        \centering
        \includegraphics[width=\linewidth]{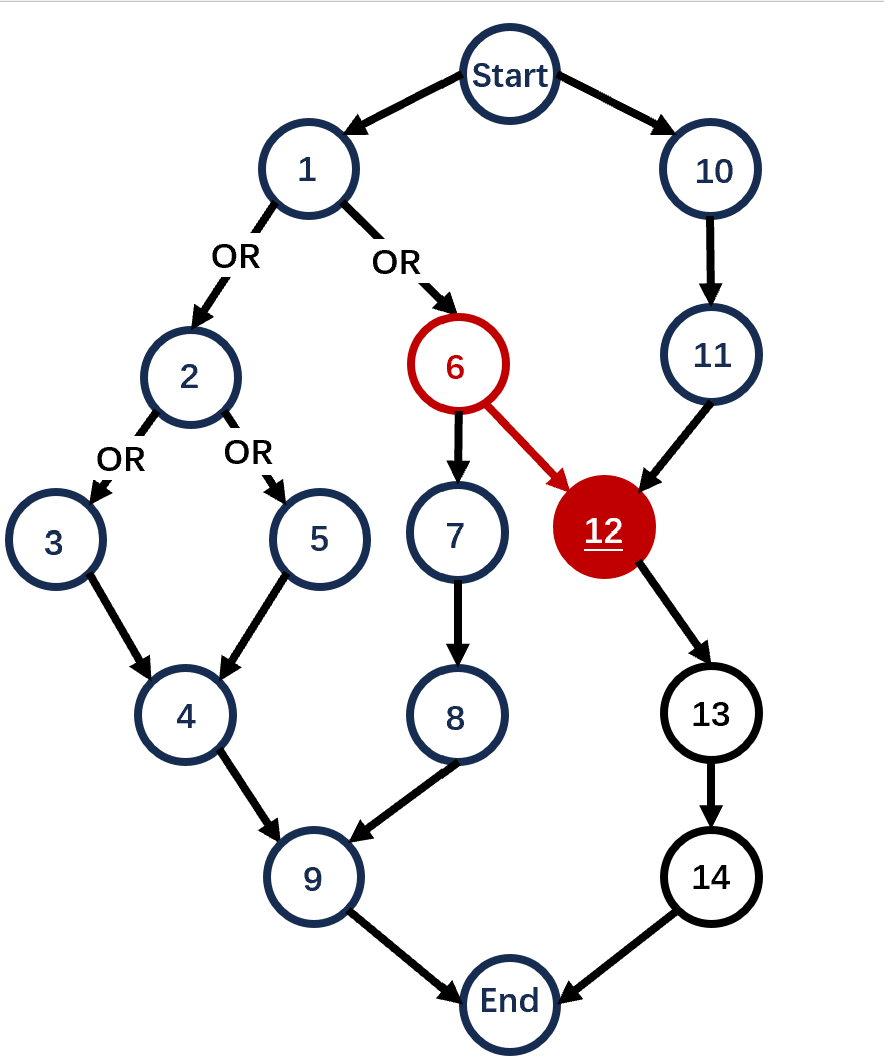}
        \caption{Sub-path to main path}
        \label{fig:sub_to_main}
    \end{subfigure}
    \caption{Comparison of path transitions.}
    \label{fig:path_transitions}
\end{figure}

In the following 2 scenarios, there are links between sub-paths and main-path such that the order of operation cannot be determined by definition.
\begin{itemize}
    \item \textbf{Link from main-path to sub-path}: It is unclear whether operation 12 can be executed as it must be completed after operation 7 and 11 while operation 12 is in an OR-link path (as shown in Figure \ref{fig:main_to_sub}).
    \item \textbf{Link from sub-path to main-path}: It's unclear that operation 12 can be done when operation 6 is not completed. (as shown in Figure \ref{fig:sub_to_main}).
\end{itemize}

To avoid the scenarios above where operations in the main path and sub-path are linked, we employ a nested approach to generate the problem incrementally. The specific steps are as follows.
\subsection*{Well-Defined IPPS Problems Generation}

    \textbf{Step 1: Construct Directed Acyclic Graph (DAG) as Main-path}
    \begin{itemize}
        \item Input: Range of the number of operations in main-path, Maximum Out-degree...
        \item Initialize the start node and end node, then establish a sequence of operation nodes between the start and end nodes under input requirements.
    \end{itemize}
    \textbf{Step 2: Add Sub-path to Main-path}
    \begin{itemize}
        \item Input: The number of OR-connector, range of the number of total operations and operations in each sub-path, Maximum number of OR-connector's out-degree
        \item Randomly select $n$ edges from the main path as initiation and integration points for the OR sub-paths, then generate \(m_i\) sub-paths with \(k_i\) nodes between originally chosen edge \(i\).  \(n, m, k\) is randomly decided under input requirements.
        \item For each added OR-path, randomly add and-path into it by selecting edges in sub-paths then generate and-path between the chosen edge. 
        \item Every time insert paths to a chosen edges, if an edge of a join-node (where two sub-path emerge) is chosen, we automatically create a supernode at the end of the new sub-path, which then connects to the original join node.
        \item Check if the number of total operations meets the constraint. If not, generate a new job. Otherwise, save the well-defined job.
    \end{itemize}
    \textbf{Step 3: Assignment of Machine and Processing Time}
    \begin{itemize}
        \item Input: The number of machines, the range of processing time and the probability that one machine is assigned to each operation.
        \item For each operation, randomly assign machines to the operation using input probabilities with processing time uniformly generated by range of processing time.
    \end{itemize}
    \textbf{Step 4: Combine Jobs to a Well-defined Problem}
    \begin{itemize}
        \item Input: The number of jobs and machines
        \item Randomly combine required jobs number of well-defined jobs using the required number of machines.
    \end{itemize}

\section{Infrastructure and Configuration}
\subsection{Infrastructure}
We use a machine with an 11th Gen Intel(R) Core(TM) i7-11800H @ 2.30GHz CPU and NVIDIA GeForce RTX 3070 GPU for DRL and greedy dispatching rules. For solving MILP by OR-tools CP-SAT Solver, we use 32 cores of an Intel Xeon Platinum 8469C at 2.60 GHz CPU with 512 GB RAM. Our core functions are implemented using PyTorch\cite{pytorch} and PyTorch Geometric(PyG)\cite{pyg}.
\subsection{Neural Networks Parameters}
Below are the essential settings for our neural networks. Note that the length of lists in the Value column for hidden\_dims or num\_heads indicates the depth.
\\

\begin{tabular}{|p{5cm}|p{1.7cm}|}
\hline
\textbf{Parameter} & \textbf{Value} \\
\hline
actor hidden\_dims (layers) & [64, 32] \\
\hline
critic hidden\_dims (layers) & [64, 32] \\
\hline
GNN pooling\_method & mean \\
\hline
GNN jumping\_knowledge & max \\
\hline
GNN model & GAT \\
\hline
GNN num\_heads (layers) & [2, 2, 2] \\
\hline
GNN hidden\_dim & 64 \\
\hline
\end{tabular}
\subsection{Environment Parameters}
We run 3 instances in a batch parallel and the size of validation set is set to 20.
\\

\begin{tabular}{|p{5cm}|p{1.7cm}|}
\hline
\textbf{Parameter} & \textbf{Value} \\
\hline
batch\_size & 3 \\
\hline
valid\_batch\_size & 20 \\
\hline
\end{tabular}

\subsection{Training Parameters}
Below are the training configurations. The learning rate (lr) and betas are fundamental parameters of the Adam optimizer. The parameter eps\_clip in PPO-clip is adjusted by adding a multiplier (clip\_multi) and an upper bound (clip\_ub) to increase it to the upper bound. It's important to note that lr, gamma, and parallel\_iter(interval of changing instances) were tuned, and we identified a range that yielded competitive results in our tuning experiment. In our final trained model, we selected values of 5e-4, 0.99, and 20 respectively.

\begin{tabular}{|p{5cm}|p{1.7cm}|}
\hline
\textbf{Parameter} & \textbf{Value} \\
\hline
lr & 5e-5 - 5e-4 \\
\hline
betas & [0.9, 0.999] \\
\hline
gamma (Discounting rewards) & 0.95 - 0.995 \\
\hline
A\_coeff (weight of policy loss) & 1 \\
\hline
clip\_ub & 0.5 \\
\hline
K\_epochs (update epoch in PPO) & 3 \\
\hline
eps\_clip & 0.3 \\
\hline
vf\_coeff (weight of value loss) & 0.5 \\
\hline
clip\_multi & 1.002 \\
\hline
entro\_coeff (weight of entropy loss) & 0.05 \\
\hline
parallel\_iter& 20 - 50 \\
\hline
save\_timestep (validate interval) & 10 \\
\hline
max\_iterations & 30000 \\
\hline
minibatch (update batchsize in PPO) & 512 \\
\hline
update\_timestep (update interval)& 5 \\
\hline
entropy\_discount & 0.997 \\
\hline
\end{tabular}
\\

\section{More Experiment Results}

\subsection{Ablation Experiments}
\begin{figure}[ht]
    \centering
    \begin{subfigure}[b]{0.45\linewidth}
        \centering
        \includegraphics[width=\linewidth]{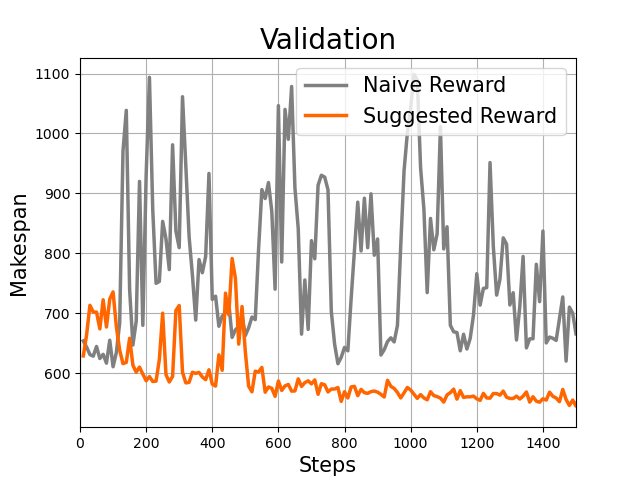}
        \caption{Reward Function}
        \label{fig:reward}
    \end{subfigure}
    \begin{subfigure}[b]{0.45\linewidth}
        \centering
        \includegraphics[width=\linewidth]{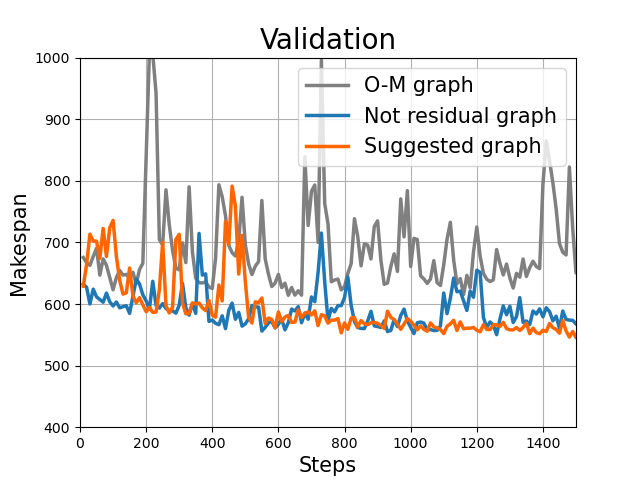}
        \caption{Graph Design}
        \label{fig:graph}
    \end{subfigure}
    \caption{Comparison of path transitions.}
    \label{fig:Alation experiment}
\end{figure}

To demonstrate the effectiveness of our reward function design along with the heterogeneous graph and updating design, we perform two ablation experiments using datasets with (5 jobs, 3 machines) as described previously. In these ablation experiments, all curves maintain the same configurations and run for 1500 training epochs, and we assess the validation curve (per 10 epochs).

Figure\ref{fig:reward} presents a comparison between the Naive Reward, which utilizes the difference in the maximum finished time of the scheduled operations, and our estimation-based Reward. Our reward function shows improved stability and competitiveness.

Figure\ref{fig:graph} presents a comparison of three distinct graph-based state representations: the \textbf{O-M graph}, which utilizes a graph with only two types of nodes — operation and machine, the \textbf{Not-residual graph} that retains nodes which are ineligible or unnecessary, and our proposed graph with four types of nodes, which exclusively keeps eligible nodes as previously detailed. The proposed graph shows enhanced competitiveness and stability. Furthermore, it achieves a faster inference speed than the \textbf{Not-residual graph} due to the reduced aggregation of neighborhood information.

\subsection{Additional Experiment in Large-scale Problem}
As shown in Table~\ref{tab:kim-gap}, the difference between the results of DRL and Greedy decreases as the size of the problem increases, indicating that our methods are particularly effective for large-scale problems. We present experimental results from two relatively large, generated datasets: 16 
\(\times\) 20 and 20 \(\times\) 25 in the main text. To further assess our efficiency, we conduct an additional experiment on 10 instances, each with 25 jobs and 20 machines. The results, shown in Table~\ref{tab:large-scale}, further confirm that our method has significant potential for large-scale problem solving.

\section{Details of Greedy Dispatching Rules}
\subsection{Selecting Operation Rules}
We utilize established greedy rules for FJSP within IPPS. To implement these FJSP rules, we initially choose a single combination for each job before application, as the rule is only executable with one specific combination.

In addition, we use a greedy rule designed for IPPS suggested by \cite{ausaf2015pgha}, which we call Muhammad.
\begin{itemize}
    \item \textbf{First In First Out (FIFO)}, which handles the earliest operation in a machine's queue to be processed.
    \item \textbf{Most Work Remaining (MWKR)}, which handles the job that has the largest amount of work left to be done.
    \item \textbf{Most Operation Number Remaining (MOR)}, which handles the job with the greatest number of remaining tasks.
\item \textbf{Muhammad}, which selects a job according to the following probability:
    \[
        P_i = \frac{C_i}{\sum_i C_i}
    \]
    where \(C_i\) is determined based on the ratio of the estimated finish time of the job \(T_m\) to the critical time \(T_c\). The definitions of \(T_m\) and \(T_c\) are as follows:

    \(T_m\) refers to the minimum machining time calculated by summing the shortest processing times for all operations in a job, and the job with the largest \(T_m\) is identified as the critical job, with its machining time defined as the critical time \(T_c\).

    The values of \(C_i\) are then assigned according to the following conditions:
    \begin{center}
    \begin{tabular}{|c|c|}
        \hline
        \textbf{Condition} & \textbf{C} \\
        \hline
        $T_m \geq 0.95T_c$ & 5 \\
        \hline
        $0.85T_c \leq T_m < 0.95T_c$ & 4 \\
        \hline
        $0.70T_c \leq T_m < 0.85T_c$ & 3 \\
        \hline
        $0.50T_c \leq T_m < 0.70T_c$ & 2 \\
        \hline
        $T_m < 0.5T_c$ & 1 \\
        \hline
    \end{tabular}
    \end{center}
    To assign priorities, each job's machining time \(T_m\) is compared with \(T_c\), and a critical value \(C_i\) is assigned based on the conditions outlined in the table.
    \end{itemize}

\subsection{Selecting Machine Rules}
For the machine selection criteria, the following rules have been considered:

\begin{itemize}
    \item \textbf{Shortest Processing Time (SPT)}, which assigns a machine with the shortest processing time for an operation.
    \item \textbf{Earliest Ending Time (EET)}, which selects the machine where an operation can start earliest.
    \item \textbf{Least Utilized Machine (LUM)}, which selects the machine with least load.
\end{itemize}

\subsection{Workflow}
The Algorithm \ref{alg:gre} displays the workflow of greedy rules. Each time an O-M pair is selected, we assign its start time to the earliest possible moment it can begin processing.
\begin{algorithm}
    \caption{Geedy Dispatching Rules}
\label{alg:gre}
\begin{algorithmic}[1]
\Require Repeat times \(\mathcal{N}\).
\State Loading instances and Initializing environment
    \For{$n = 1, 2, \dots, \mathcal{N}$} 
        \While{$s_t$ is not terminal}
            \State Use Greedy Rules to select a O-M pair.
            \State Schedule the pair at the earliest possible time.
            \State Receive reward $r_t$ and enter the next state $s_{t+1}$.
            \State Update environment.
        \EndWhile
    \EndFor
\State Choose the solution with min makespan.

\State \textbf{return}
\end{algorithmic}
\end{algorithm}

\section{Experiment Details in Open Benchmarks}

\subsection{Algorithms}
In this section, we present a detailed training and testing process.
\subsubsection{Training Phase}
For training, we use the PPO-clip algorithm. The pseudo-code of PPO algorithm is as follows. 

\begin{algorithm}
\caption{Training procedure of PPO-clip}
\label{alg:ppo_training}
\begin{algorithmic}[1]

\Require Total training epoch is \(\mathcal{T}\), the interval of changing instances is \(\mathcal{P}\) and the interval of updating is \(\mathcal{U}\)
\State Initializing trainable parameters \(\omega\) of HGAT, policy network \(\textbf{Actor}_\theta\) and Value \(\textbf{Critic}_\phi\) network.
\State Initializing Data-Generator or preparing set of IPPS instances.
\For{$\text{iter} = 1, 2, \dots, \mathcal{T}$}
    \State Loading data to a batch and Initializing environment
    \For{$b = 1, 2, \dots, \mathcal{P}$} 
        \While{$s_t$ is not terminal}
            \State Extract embeddings of O-M pairs.
            \State Sample $a_t \sim \textbf{Actor}_\theta(\cdot | s_t)$
            \State Receive reward $r_t$ and enter next state $s_{t+1}$
            \State Update environment and memory.
        \EndWhile
        \If{b mod \(\mathcal{U}\) = 0}
            \State Flatten states in memory.
            \State Calculate PPO loss.
            \State Update network parameters.
        \EndIf
    \EndFor

\EndFor
\State \textbf{return}
\end{algorithmic}
\end{algorithm}

\subsubsection{Inference Phase}
DRL-Sample (DRL-S) and DRL-Greedy (DRL-G) are two methods used during the inference phase. Specifically, DRL-G selects actions based on the highest probability, whereas DRL-S selects actions according to the probability distribution determined by the actor. The pseudo-code of these methods are as follows:
\begin{algorithm}
    \caption{DRL-G}
\label{alg:drl-g}
\begin{algorithmic}[1]
\Require Repeat times \(\mathcal{N}\).
\State Load the models that have already been trained.
\State Loading instances and Initializing environment
    \For{$n = 1, 2, \dots, \mathcal{N}$} 
        \While{$s_t$ is not terminal}
            \State Extract embeddings of O-M pairs.
            \State Select $a_t  = \textbf{argmax}(\textbf{Actor}_\theta(\cdot | s_t))$
            \State Get the next state $s_{t+1}$
            \State Update environment.
        \EndWhile
    \EndFor
\State Choose the solution with min makespan.

\State \textbf{return}
\end{algorithmic}
\end{algorithm}
\begin{algorithm}
    \caption{DRL-S}
\label{alg:drl-s}
\begin{algorithmic}[1]
\Require Repeat times \(\mathcal{N}\).
\State Load the models that have already been trained.
\State Loading instances and Initializing environment
    \For{$n = 1, 2, \dots, \mathcal{N}$} 
        \While{$s_t$ is not terminal}
            \State Extract embeddings of O-M pairs.
            \State Sample $a_t \sim \textbf{Actor}_\theta(\cdot | s_t)$
            \State Get next state $s_{t+1}$
            \State Update environment.
        \EndWhile
    \EndFor
\State Choose the solution with min makespan.

\State \textbf{return}
\end{algorithmic}
\end{algorithm}

\subsection{Results}

 The complete results of the Kim dataset are shown in Table~\ref{tab:kim-g} and Table~\ref{tab:kim-drl}.
\begin{table}[ht]
\centering
\resizebox{0.7\textwidth}{!}{
\begin{tabular}{ccccccc}
\toprule
\textbf{Size} & \multicolumn{3}{c}{\textbf{DRL-G}} & \multicolumn{3}{c}{\textbf{DRL-S}} \\
\cmidrule(lr){2-4} \cmidrule(lr){5-7}
 & \textbf{\(\mathbf{4\times5}\)} & \textbf{\(\mathbf{5\times3}\)} & \textbf{\(\mathbf{6\times5}\)} & \textbf{\(\mathbf{4\times5}\)} & \textbf{\(\mathbf{5\times3}\)} & \textbf{\(\mathbf{6\times5}\)} \\
\midrule
\(6\times 15\)  & 6.04\%  & 7.23\%  & 7.19\%  & 0.13\%   & -0.18\% & 0.63\%  \\
\(9\times 15\)  & -0.08\% & 2.17\%  & -1.05\% & -6.25\%  & -5.75\% & -6.29\% \\
\(12\times 15\) & -1.84\% & 0.09\%  & -3.60\% & -9.80\%  & -10.15\%& -10.70\%\\
\(15\times 15\) & -6.29\% & -9.19\% & -5.79\% & -14.62\% & -12.80\%& -14.42\%\\
\(18\times 15\) & -10.71\%& -4.54\% & -6.72\% & -10.16\% & -10.34\%& -11.07\%\\
\bottomrule
\end{tabular}
}

\caption{Average Gap Across Different Problem Sizes in the Kim dataset. The gap is relative to the best greedy rule for each problem.}
\label{tab:kim-gap}
\end{table}

\begin{table*}[ht]
\centering
\resizebox{\textwidth}{6cm}{
\begin{tabularx}{\linewidth}{lXXXXXXXXXXXXX}
\toprule
\multirow{2}{*}{Greedy rules} & \multicolumn{3}{c}{MWKR} & \multicolumn{3}{c}{MOR} & \multicolumn{3}{c}{FIFO} & \multicolumn{3}{c}{Muhammad} & \multirow{2}{*}{BCI} \\
\cmidrule(lr){2-4} \cmidrule(lr){5-7} \cmidrule(lr){8-10} \cmidrule(lr){11-13}
 & SPT & EFT & LUM & SPT & EFT & LUM & SPT & EFT & LUM & SPT & EFT & LUM & \\
\midrule
problem01 & 452 & 523 & 543 & 427 & 532 & 518 & 475 & 544 & 538 & 427 & 472 & 494 & 427 \\
problem02 & 422 & 475 & 460 & 387 & 392 & 393 & 438 & 435 & 431 & 398 & 376 & 387 & 376 \\
problem03 & 417 & 473 & 447 & 385 & 422 & 420 & 469 & 491 & 504 & 374 & 388 & 378 & 374 \\
problem04 & 369 & 381 & 359 & 319 & 329 & 331 & 382 & 361 & 371 & 311 & 330 & 336 & 311 \\
problem05 & 419 & 555 & 455 & 386 & 388 & 392 & 423 & 388 & 414 & 375 & 364 & 369 & 364 \\
problem06 & 571 & 537 & 509 & 558 & 508 & 491 & 504 & 536 & 572 & 526 & 505 & 489 & 489 \\
problem07 & 447 & 500 & 571 & 434 & 401 & 434 & 443 & 456 & 483 & 425 & 392 & 388 & 388 \\
problem08 & 536 & 448 & 429 & 432 & 393 & 401 & 443 & 403 & 414 & 427 & 378 & 395 & 378 \\
problem09 & 494 & 523 & 505 & 544 & 511 & 535 & 427 & 533 & 518 & 429 & 486 & 474 & 427 \\
problem10 & 540 & 537 & 521 & 492 & 510 & 489 & 567 & 553 & 564 & 505 & 506 & 509 & 489 \\
problem11 & 536 & 513 & 495 & 463 & 420 & 433 & 502 & 470 & 521 & 447 & 442 & 446 & 420 \\
problem12 & 446 & 477 & 465 & 402 & 353 & 399 & 505 & 398 & 445 & 391 & 377 & 400 & 353 \\
problem13 & 653 & 610 & 596 & 563 & 500 & 503 & 641 & 566 & 618 & 564 & 521 & 503 & 500 \\
problem14 & 555 & 614 & 617 & 517 & 441 & 487 & 496 & 521 & 540 & 500 & 456 & 461 & 441 \\
problem15 & 494 & 571 & 505 & 555 & 499 & 511 & 490 & 514 & 537 & 483 & 496 & 485 & 483 \\
problem16 & 581 & 583 & 586 & 533 & 490 & 507 & 601 & 585 & 627 & 570 & 529 & 547 & 490 \\
problem17 & 614 & 618 & 551 & 543 & 452 & 490 & 551 & 547 & 517 & 540 & 518 & 517 & 452 \\
problem18 & 503 & 533 & 512 & 488 & 391 & 429 & 524 & 476 & 476 & 454 & 445 & 457 & 391 \\
problem19 & 677 & 663 & 654 & 604 & 513 & 540 & 655 & 606 & 637 & 596 & 559 & 566 & 513 \\
problem20 & 618 & 590 & 571 & 539 & 475 & 488 & 579 & 552 & 571 & 577 & 532 & 543 & 475 \\
problem21 & 662 & 669 & 621 & 614 & 508 & 512 & 565 & 574 & 648 & 563 & 506 & 538 & 506 \\
problem22 & 709 & 721 & 718 & 673 & 535 & 557 & 640 & 622 & 693 & 630 & 587 & 601 & 535 \\
problem23 & 642 & 662 & 628 & 579 & 487 & 517 & 621 & 639 & 641 & 562 & 566 & 570 & 487 \\
problem24 & 770 & 756 & 740 & 709 & 551 & 607 & 754 & 714 & 716 & 696 & 636 & 648 & 551 \\
\midrule
Average & 546.95 & 563.83 & 544.08 & 506.08 & 458.38 & 474.33 & 528.96 & 520.17 & 541.50 & 490.47 & 473.63 & 479.21 & 458.38 \\
\bottomrule
\end{tabularx}
}

\caption{Greedy Result for Kim.}
\label{tab:kim-g}
\end{table*}

\begin{table*}[ht]
\centering
\resizebox{\textwidth}{5.5cm}{
\begin{tabularx}{\linewidth}{lp{1cm}p{1cm}p{1cm}p{1cm}p{1cm}p{1cm}p{1cm}p{1cm}XX}
\toprule
\textbf{Problem} & \multicolumn{3}{c}{\textbf{DRL-G}} & \multicolumn{3}{c}{\textbf{DRL-S}} & \multicolumn{2}{c}{\textbf{Greedy}} &\multirow{2}{*}{\textbf{OR-tool}}&\multirow{2}{*}{\textbf{LB}{*}} \\
\cmidrule(lr){2-4} \cmidrule(lr){5-7} \cmidrule(lr){8-9}
& \textbf{4 \(\times\) 5} & \textbf{5\(\times\)3} & \textbf{6\(\times\)5} & \textbf{4\(\times\)5} & \textbf{5\(\times\)3} & \textbf{6\(\times\)5} & \textbf{Best One}{\dag} & \textbf{Best Choice}{\ddag} & &\\
\midrule
problem01 & 523 & 496 & 486 & 462 & 468 & 484 & 532 & \textbf{427} & 427 & 427 \\
problem02 & 377 & 413 & 421 & \textbf{360} & 363 & 364 & 392 & 376 & 343 & 343  \\
problem03 & 372 & 375 & 390 & 372 & \textbf{363} & 367 & 422 & 374 & 346 & 344  \\
problem04 & 334 & 364 & 325 & 318 & 312 & 315 & 329 & \textbf{311} & 306 & 306  \\
problem05 & 343 & 338 & 345 & 332 & \textbf{330} & 331 & 388 & 364 & 318 & 318  \\
problem06 & 528 & 496 & 533 & 483 & \textbf{474} & 487 & 508 & 489 & 427 & 427  \\
problem07 & 396 & 403 & 392 & \textbf{374} & 383 & 385 & 401 & 388 & 372 & 372  \\
problem08 & 397 & 407 & 404 & 362 & 360 & \textbf{356} & 393 & 378 & 343 & 343  \\
problem09 & 493 & 495 & 506 & 484 & 484 & 479 & 511 & \textbf{427} & 427 & 427  \\
problem10 & 517 & 544 & 500 & \textbf{474} & 486 & \textbf{474} & 510 & 489 & 427 & 427  \\
problem11 & 371 & 441 & 404 & \textbf{367} & \textbf{367} & 368 & 420 & 420 & 344 & 344  \\
problem12 & 357 & 338 & 367 & 338 & 334 & \textbf{332} & 353 & 353 & 318 & 318  \\
problem13 & 525 & 498 & 492 & \textbf{478} & 479 & 485 & 500 & 500 & 427 & 427  \\
problem14 & 407 & 429 & 405 & \textbf{383} & 385 & 384 & 441 & 441 & 372 & 372  \\
problem15 & 517 & 503 & 488 & 483 & 488 & \textbf{481} & 499 & 483 & 427 & 427  \\
problem16 & 510 & 518 & 505 & 483 & 478 & \textbf{475} & 490 & 490 & 427 & 427  \\
problem17 & 392 & 419 & 405 & 377 & \textbf{372} & 373 & 452 & 452 & 382 & 344  \\
problem18 & 411 & 390 & 385 & 350 & 353 & \textbf{344} & 391 & 391 & 322 & 318  \\
problem19 & 504 & 557 & 502 & \textbf{472} & 473 & 477 & 513 & 513 & 427 & 427  \\
problem20 & 439 & 424 & 434 & 399 & 389 & \textbf{387} & 475 & 475 & 385 & 372  \\
problem21 & 518 & 529 & 496 & \textbf{474} & 480 & 475 & 508 & 506 & 427 & 427  \\
problem22 & 538 & 496 & 506 & \textbf{473} & 487 & \textbf{473} & 535 & 535 & 427 & 427  \\
problem23 & 423 & 433 & 457 & \textbf{401} & 406 & 403 & 487 & 487 & 534 & 372  \\
problem24 & 492 & 526 & 514 & 495 & 494 & \textbf{490} & 551 & 551 & 455 & 427  \\
\bottomrule
\end{tabularx}
}
\caption{DRL Result For Kim. \dag Best One refers to the result of greedy rule that has a minimum average makespan while \ddag Best Choice is the best greedy result for each problem. *LB refers to lower bound of Kim dataset demonstrated in \cite{Liu2021}}
\label{tab:kim-drl}
\end{table*}
\begin{table*}[ht]
\centering
\begin{tabularx}{\textwidth}{lXXXXXXXXXXXXX}
\toprule
\multirow{2}{*}\textbf{Problem} & \multicolumn{3}{c}{\textbf{DRL-S}} & \multicolumn{3}{c}{\textbf{DRL-G}} & \multicolumn{2}{c}{\textbf{Greedy}} & \multirow{2}{*}{\textbf{OR-tool}} \\
\cmidrule(lr){2-4} \cmidrule(lr){5-7} \cmidrule(lr){8-9}
 & \textbf{4 \(\times\) 5} & \textbf{5\(\times\)3} & \textbf{6\(\times\)5} & \textbf{4\(\times\)5} & \textbf{5\(\times\)3} & \textbf{6\(\times\)5} & \textbf{Best One}{\dag} & \textbf{Best Choice}{\ddag} &  \\
\midrule
problem01  & \textbf{412} & 424 & 415 & 427 & 446 & 438 & 480 & 480 & 621 \\
problem02  & 417 & \textbf{412} & 423 & 424 & 438 & 781 & 473 & 473 & 1659 \\
problem03  & 421 & 419 & \textbf{415} & 496 & 438 & 499 & 483 & 483 & 692 \\
problem04  & 418 & \textbf{413} & 423 & 459 & 444 & 452 & 479 & 479 & 1045 \\
problem05  & 417 & 421 & \textbf{419} & 537 & 420 & 434 & 480 & 480 & 1134 \\
problem06  & 421 & 417 & \textbf{414} & 448 & 449 & 430 & 479 & 479 & 887 \\
problem07  & \textbf{417} & 419 & 416 & 429 & 421 & 429 & 472 & 472 & 481 \\
problem08  & 420 & 420 & \textbf{419} & 465 & 440 & 442 & 476 & 476 & 1152 \\
problem09  & 422 & 419 & \textbf{418} & 438 & 432 & 452 & 465 & 465 & 884 \\
problem10 & 419 & \textbf{410} & 425 & 445 & 439 & 455 & 473 & 473 & 469 \\
\midrule
Average Gap & -46.32\% & \textbf{-46.43\%} & -46.28\% & -41.82\% & -44.01\% & -40.46\% & -38.93\% & -38.93\% & - \\
\bottomrule
\end{tabularx}
\caption{Performance in problems of size 25 $\times$ 20. \dag Best One refers to the result of greedy rule that has a minimum average makespan while \ddag Best Choice is the best greedy result for each problem. }
\label{tab:large-scale}
\end{table*}

In Table~\ref{tab:kim-gap}, we compare the gap relative to the best greedy result for each problem across different sizes. It is evident that the gap between the results of DRL and Greedy decreases as the problem size increases. This suggests that our method demonstrates significant potential in solving large-scale problems.

\end{appendices}
\end{document}